\documentclass[10pt,openbib]{article}

\usepackage[T1]{fontenc}
\usepackage{amsfonts,amsthm}
\usepackage{amsmath}
\usepackage{amssymb}

\newcommand{\R}{\mathbb{R}}

\newcommand{\Z}{\mathbb{Z}}
\newcommand{\Q}{\mathbb{Q}}

\newcommand{\onehalf}{\mbox{${\scriptstyle \frac{1}{2}\, }$}}
\newcommand{\dee}{\mathop{\! \, \rm d \!}\nolimits}
\newcommand{\comp}{\, \raisebox{2pt}{$\scriptstyle\circ \, $}}
\newcommand{\setrule}{\, \rule[-4pt]{.5pt}{13pt}\, }
\newcommand{\rowspace}{\rule{0pt}{16pt}}
\newcommand{\spann}{\mathop{\rm span}\nolimits}
\newcommand{\lefthook}{\mbox{$\, \rule{8pt}{.5pt}\rule{.5pt}{6pt}\, \, $}}

\newcommand{\dbyds}{\mbox{${\displaystyle \frac{\dee}{\dee t}}
\rule[-10pt]{.5pt}{25pt} \raisebox{-10pt}{$\, {\scriptstyle s=0}$}$}}
\newcommand{\smalldbydt}{\mbox{ {\footnotesize $\frac{\dee}{\dee t}
\rule[-8pt]{.35pt}{15pt} \raisebox{-8pt}{$\, {\scriptstyle t=0}$}$ }}}

\allowdisplaybreaks

\begin{document}

\title{\textbf{Vector fields and differential forms on the} \\ \textbf{orbit space of a proper action.} } 
\author{Larry Bates, Richard Cushman, and J\k{e}drzej \'{S}niatycki\footnotemark }
\addtocounter{footnote}{1}
\footnotetext{email. bates@ucalgary.ca, r.h.cushman@gmail.com, sniatycki@gmail.com}
\date{}
\maketitle
\addtocounter{footnote}{2}
\footnotetext{printed: \today}

\begin{abstract}
In this paper we study differential forms and vector fields on the \linebreak 
orbit space of a proper action of a Lie group  on a smooth manifold, defining them as multilinear maps on the generators of infinitesimal diffeomorphisms, respectively. This yields an intrinsic view of vector fields and differential forms on the  orbit space. 
\end{abstract}

\section{Introduction}

This paper is a complete version of \cite{bates-cushman-sniatycki21} 
which includes the correction to Lemma 3.1, whose proof contained 
an error, which was pointed out by Prof. G. Schwarz. \medskip 

\vspace{-.15in}This paper is part of a series of papers devoted to the study of the 
geometry of singular spaces in terms of the theory of differential spaces, 
which were introduced by Sikorski \cite{sikorski67}, see also \cite{sikorski72}. 
In this theory, geometric information about a space $S$ is encoded in a 
ring $C^{\infty}(S)$ of real valued functions, which are deemed to be smooth. In 
particular, we are concerned with the class of subcartesian spaces introduced 
by Aronszajn \cite{aronszajn67}. A Hausdorff differential space $S$ is \emph{subcartesian} 
if every point $x$ of $S$ has a neighborhood $U$ that is diffeomorphic to 
a subset $V$ of a Euclidean (Cartesian) space ${\R }^n$. The restriction of 
$C^{\infty}(S)$ to $U$ is isomorphic to the restriction of $C^{\infty}({\R }^n)$ to 
$V$, see \cite{sniatycki}. \medskip 

Palais \cite{palais} introduced the notion of a slice for an action of a 
not necessarily compact Lie group $G$ on a manifold $M$. Since then, the structure of the space 
$M/G$ of orbits of a proper action of $G$ on $M$ has been investigated by 
many mathematicians. In \cite{duistermaat} Duistermaat showed that $M/G$ is 
a subcartesian differential space with differential structure $C^{\infty}(M/G)$ 
consisting of push forwards of smooth $G$ invariant functions on $M$ by the 
$G$ orbit mapping $\pi : M \rightarrow M/G$. \medskip 

\vspace{-.15in}On a smooth manifold $M$, there are two equivalent definitions of a vector field, namely,  as a derivation of $C^{\infty}(M)$, or as a generator of a local one parameter local group of diffeomorphisms of $M$. Choosing one, and proving the other is a matter of preference. 
On a subcartesian differential space $S$, which is not in general a manifold, these notions differ. We use the term \textit{vector field} on $S$ for a generator of a local one parameter group of local diffeomorphisms and denote the class of all vector fields on $S$ by $\mathfrak{X}(S)$. A key reason for the choice made in this paper is the special case of the orbit space of a proper action. The class of derivations of $C^{\infty}(S)$ is, in general, larger than the class 
$\mathfrak{X}(S)$. For $S = M/G$ we show that a derivation $Y$ of 
$C^{\infty}(M/G)$ is in 
$\mathfrak{X}(S)$, that is, is a vector field on $M/G$, if and only if there exists 
a $G$ invariant vector field $X$ on $M$ such that $Y$ is $\pi $ related to $X$, that is, $Y \comp \pi = T\pi \comp X$, where $\pi : M \rightarrow M/G$ is the $G$ orbit map. 
\medskip %

In the literature there has been extensive discussion about the notion of 
a differential form on a singular space, see Smith \cite{smith}, Marshall 
\cite{marshall75}, and Sjamaar \cite{sjamaar}. Here, in our search for an intrinsic notion of a differential form, we have been led to see them 
as multilinear maps on vector fields. In the case of a $1$-form $\theta $ on $M/G$ with $\theta $ a linear mapping  
\begin{displaymath}
\theta : \mathfrak{X}(M/G) \rightarrow C^{\infty}(M/G): Y \mapsto 
\langle \theta | Y \rangle 
\end{displaymath}%
over the ring $C^{\infty}(M/G)$ of smooth functions on $M/G$, which is to say  
$\langle \theta | \overline{f}Y \rangle = 
\overline{f} \langle \theta | Y \rangle $ for every $\overline{f} \in C^{\infty}(M/G)$. 
With this definition we show that every differential $1$-form on $M/G$ pulls back under the orbit map $\pi $ to a semi-basic $G$ invariant $1$-form on $M$. Furthermore every $G$ invariant semi-basic differential $1$-form on $M$ is the pull back by $\pi $ of a differential $1$-form on $M/G$. \medskip 

We define a differential exterior algebra of differential forms on the orbit space, which satisfies a version of de Rham's theorem. Our version is larger than Smith's as it includes forms that are not Smith forms, see section 6. It also handles singular orbit spaces of a proper action of a Lie group on a smooth manifold. The Lie group need not be compact and the orbit space need not be smooth, both of which Koszul hypothesized in \cite{koszul}. \medskip

We now give a section by section description of the contents of this paper. \medskip 

Section 2 deals with basic properties of a proper action of a Lie group $G$ on 
a smooth manifold $M$ and the differential structure of the orbit space $M/G$. 
We introduce the reader to the theory of subcartesian
differential spaces in the context of the orbit space $M/G$. The differential geometry of $M$ is described in terms of its smooth structure given by the ring 
$C^{\infty}(M)$ of smooth functions on $M$. The differential geometry of the 
orbit space $M/G$, which may have singularities, is similarly described in terms of the ring $C^{\infty }(M/G)$ of smooth functions on $M/G$, which is isomorphic to the ring $C^{\infty }(M)^{G}$ of smooth, $G$-invariant functions on $M$. 
Since a proper action has an invariant Riemannian metric, several results are proved using properties of the geodesics of the metric.  Also certain objects are shown to be smooth submanifolds. \medskip 

In section 3 we study vector fields on the orbit space $M/G$. In the case of the manifold $M$, derivations of the ring $C^{\infty }(M)$ are vector fields on $M$, 
and they generate a local one parameter group of local diffeomorphisms of 
$M$. In the case of the ring $C^{\infty }(M/G)$, not all derivations of 
$C^{\infty }(M/G)$ generate local one parameter groups of local diffeomorphisms of $M/G$. The derivations of $C^{\infty }(M/G),$ which generate local one parameter groups of local diffeomorphisms of $M/G$. This is the key idea of this paper.  We establish that every vector field on the quotient $M/G$ is covered by a $G$-invariant vector field on $M$. It is well known that the space $M/G$ is stratified, see \cite{duistermaat,thom,whitney}. We show that every vector field on $M/G$ defines a vector field tangent to each stratum of $M/G$.\medskip 

In section 4, we define differential $1$-forms on the orbit space $M/G$ as
linear mappings on the the space $\mathfrak{X}(M/G)$ of smooth vector fields
on $M/G$. The most important consequence of this definition relates to pulling back $1$-forms from $M/G$ to $M$.  In particular, our notion of a differential 
$1$-form is intrinsic. \medskip 

In section 5 to prove a version of de Rham's theorem we enlarge  
the algebra of differential $1$-forms to $k$-forms with an 
exterior derivative operator.  The key techinical point is that everything is developed in terms of the Lie derivative of vector fields.  Almost all of this section looks the same as that on manifolds. \medskip 

In section 6 we give all the details of the simplest nontrivial example. 
This example reveals that differential forms in our sense are not the same as 
those of Smith \cite{smith}.

\section{Basic properties}

This section gives some of the basic properties of smooth vector fields on the 
orbit space of a proper action of a Lie group on a smooth manifold. \medskip

Let $M$ be a connected smooth manifold with a proper action 
\begin{equation}
\Phi : G \times M \rightarrow M: (g,m) \mapsto {\Phi}_g(m) = g \cdot m
\label{eq-s0zero}
\end{equation}
of a Lie group $G$ on $M$, and let 
\begin{equation}
\pi : M \rightarrow M/G: m \mapsto \overline{m} = G \cdot m = \{ {\Phi }_g(m) \in M \setrule \, 
g \in G \}
\label{eq-s0zeronw}
\end{equation}
be the orbit map of the $G$ action $\Phi $. \medskip 

Let $C^{\infty}(M)^G$ be the algebra 
of smooth $G$ invariant functions on $M$ and let $C^{\infty}(M/G)$ be the 
algebra of functions $\overline{f}$ on $M/G$ such that $f = {\pi }^{\ast }\overline{f} = 
\overline{f} \comp \pi $ lies in $C^{\infty}(M)^G$. The map 
\begin{displaymath}
{\pi }^{\ast }: C^{\infty}(M/G) \rightarrow C^{\infty}(M)^G: \overline{f} \mapsto \overline{f} \comp \pi 
\end{displaymath}
is a bijective algebra isomorphism, whose inverse is 
\begin{displaymath}
{\pi }_{\ast }: C^{\infty}(M)^G \rightarrow C^{\infty}(M/G): f \mapsto \overline{f}. 
\end{displaymath}
  
\noindent \textbf{Proposition 1.} \textit{The orbit space $M/G$ with the differential 
structure $C^{\infty}(M/G)$ is a locally closed subcartesian differential space.} \medskip  

\noindent \textbf{Proof.} See corollary 4.11 of Duistermaat \cite{duistermaat} and 
page 72 of \cite{sniatycki}. \hfill $\square $ \medskip 

Let $X$ be a smooth vector field on a manifold $M$. $X$ gives rise to a map 
\begin{displaymath}
L_X: C^{\infty}(M) \rightarrow C^{\infty}(M): f \mapsto L_Xf = X(f), 
\end{displaymath}
called the \emph{derivation} associated to $X$. If we want to emphasize this action of vector fields on $M$, we say that they form the space $\mathrm{Der}\, C^{\infty}(M)$ of derivations of $C^{\infty}(M)$. If we want to emphasize that $X$ generates a local one parameter group of local diffeomorphisms of $M$, we say that $X$ is a 
\emph{vector field} on $M$ and write $\mathfrak{X}(M)$ for the set of vector fields on $M$. For each smooth manifold $M$ we have $\mathfrak{X}(M) = \mathrm{Der}\, C^{\infty}(M)$. However, these notions need not coincide for a subcartesian differential space. \medskip 

Let $(S, C^{\infty}(S))$ be a differential space with $X$ a derivation of 
$C^{\infty}(S)$. Let $I_x \subseteq \R \rightarrow S: t \mapsto {\varphi }^X_t(x)$ be a maximal integral curve of $X$, which starts at $x$. Here $I_x$ is an interval 
containing $0$. If $t$, $s$, and $t+s$ lie in $I_x$, and if 
$s\in I_{{\varphi }^X_t(x)}$ and $t \in I_{{\varphi }^X_s(x)}$, then 
\begin{displaymath}
{\varphi }^X_{t+s}(x) = {\varphi }^X_s({\varphi}^X_t(x)) = 
{\varphi }^X_t({\varphi }^X_s(x)).
\end{displaymath}
The map ${\varphi }^X_t$ may fail to be a local diffeomorphism of the differential 
space $S$, see example 3.2.7 in \cite[p.37]{sniatycki}.  
A \emph{vector field} on a subcartesian differential space 
$S$ is a derivation $X$ of $C^{\infty}(S)$ such that for every $x \in S$ there is 
an open neighborhood $U$ of $x$ and $\varepsilon >0$ such that 
for every $t \in (-\varepsilon , \varepsilon )$ the map ${\varphi }^X_t$ is defined on 
$U$ and its restriction to $U$ is a diffeomorphism from $U$ onto an open 
subset of $S$. In other words, the derivation $X$ is a vector field on $S$ 
if $t \mapsto {\varphi }^X_t$ is a local one parameter group of local 
diffeomorphisms of $S$. \medskip 

\noindent \textbf{Example 1.} Consider $\Q  \subseteq \R $ with the structure of 
a differential subspace of $\R $. Let $\iota : \Q  \rightarrow \R $ be the 
inclusion mapping. The differential structure $C^{\infty}(\Q )$ of $\Q $ consists of 
${\iota }^{\ast}f$, which is the restriction of a smooth function 
$f$ on $\R $ to $\Q $. Let 
$X(x_1) = a_1(x_1) \frac{\partial }{\partial x_1}$ be a vector field on $\R $. Then for every 
$f \in C^{\infty}(\R )$ and every $x_1 \in \R $ the function 
$x_1 \mapsto X(f)(x_1) = a_1(x_1) \frac{\partial f}{\partial x_1}(x_1)$ is smooth. Restricting to points $x_1$ in $\Q $ we get ${\iota }^{\ast }(X(f)) = {\iota }^{\ast }(a_1) \, {\iota }^{\ast }\big( \frac{\partial f}{\partial x_1} \big) $. We now show that we can obtain ${\iota }^{\ast }\big( \frac{\partial f}{\partial x_1} \big)$ by 
operations on ${\Q }^2$. Let $x^0_1 \in \Q $ and let $\{ (x_1)_n \} $ be a 
sequence of points in $\Q $, which converges to $x^0_1$. Then 
\begin{displaymath}
\lim_{n \rightarrow \infty} \frac{\partial f}{\partial x_1}(x^0_1+(x_1)_n) = 
\frac{\partial f}{\partial x_1}(x^0_1). 
\end{displaymath}
Thus we have shown that 
${\iota }^{\ast }(X(f)) = X_{|\Q}({\iota }^{\ast }f)$ for every $f \in C^{\infty}(\R )$. 
In other words, the restriction $X_{|\Q }$ of the vector field $X$ to $\Q $ is 
a derivation of $C^{\infty}(\Q)$. Thus $\mathrm{Der}\, C^{\infty}(\Q ) = 
\{ X_{|\Q} \setrule \, X \in \mathfrak{X}(\R ) \} $. However, no two distinct points 
of $\Q $ can be joined by a smooth curve. Hence only the derivation of $C^{\infty}(\Q )$ 
that is identically $0$ on $\Q$ admits integral curves, that is, $\mathfrak{X}(\Q ) = 
\{ 0 \} $. \hfill $\square $ \medskip 

Let $\mathfrak{X}(M)^G$ be the set of smooth $G$ invariant vector fields on $M$, 
that is,  
\begin{displaymath}
\mathfrak{X}(M)^G = \{ X \in \mathfrak{X}(M) \setrule \, T_m{\Phi }_g \, X(m)  = 
X ({\Phi }_g (m)) \, \, \mbox{for every $(g,m) \in G \times M$} \} . 
\end{displaymath}
Since $\mathfrak{X}(M) = \mathrm{Der}\, C^{\infty}(M)$, we have $\mathfrak{X}(M)^G = 
(\mathrm{Der}\, C^{\infty}(M))^G$. Also we may consider the space 
$\mathrm{Der}\,  C^{\infty}(M)^G$ of derivations of $C^{\infty}(M)^G$. Clearly, we have  $(\mathrm{Der}\, C^{\infty}(M))^G \subseteq \mathrm{Der}\, C^{\infty}(M)^G$. For $X \in \mathfrak{X}(M)^G$ let 
\begin{equation}
\varphi : D \subseteq \R \times M \rightarrow M: (t,m) \mapsto {\varphi }_t(m)
\label{eq-s0zerostar}
\end{equation}
be the \emph{local flow} of $X$, that is, $\varphi $ is a differentiable mapping such that 
\begin{displaymath}
\frac{\dee {\varphi }_t}{\dee t}(m) = X\big( {\varphi }_t(m) \big), \, \, 
\mbox{for all $(t,m) \in D$.}
\end{displaymath}
Here $D$ is a \emph{domain}, that is, $D$ is the largest (in the sense of containment) open subset of $\R \times M$ such that for each $m \in M$ the set $\{ t \in {\R } \setrule \, (t,m) \in D\} $ is an open interval containing $0$. Moreover, ${\varphi }(0,m) = m$ for every $m \in M$ and if $(t,m) \in D$, 
$(s, {\varphi }_t(m)) \in D$, and $(t+s,m) \in D$, then $ {\varphi }_{s+t}(m) = 
{\varphi }_s\big( {\varphi }_t(m) \big) $. Since $X \in \mathfrak{X}(M)^G$,  
\begin{equation}
({\Phi }_g \comp {\varphi }_t)(m) = {\varphi }_t \big( {\Phi }_g(m) \big), \, \, 
\mbox{for all $\big( g, (t,m) \big) \in G \times D$.}
\label{eq-one}
\end{equation}
Thus $(t, {\Phi }_g(m)) \in D$ for all $g \in G$, if $(t,m) \in D$. \medskip

\noindent \textbf{Proposition 2.} \textit{Let $X \in \big( \mathrm{Der}\, C^{\infty}(M) \big)^G$. Then $X$ induces a derivation of $C^{\infty}(M/G)$ defined by 
\begin{displaymath}
Y: C^{\infty}(M/G) \rightarrow C^{\infty}(M/G): \overline{f} \mapsto 
{\pi }_{\ast }\big( X ({\pi }^{\ast }\overline{f}) \big) . 
\end{displaymath}
This leads to the module homomorphism
\begin{equation}
(\mathrm{Der}\, C^{\infty}(M))^G \rightarrow \mathrm{Der}\, C^{\infty}(M/G): 
X \mapsto Y = {\pi }_{\ast } \comp X \comp {\pi }^{\ast }. 
\label{eq-s2onestarnw}
\end{equation} }

\noindent \textbf{Proof.} Let $Y = {\pi }_{\ast } \comp X \comp {\pi }^{\ast }$ and 
$\overline{f}$, $\overline{h} \in C^{\infty}(M/G)$. Then 
\begin{align*}
Y(\overline{h}\, \overline{f}) & = {\pi }_{\ast }\big( X({\pi }^{\ast }(\overline{h}\, \overline{f})) \big) 
= {\pi }_{\ast }\big( ({\pi }^{\ast }\overline{h}) X({\pi }^{\ast }\overline{f})  
+ ({\pi }^{\ast }\overline{f}) X({\pi }^{\ast }\overline{h}) \big) \\
& = \big( ({\pi }_{\ast } \comp {\pi }^{\ast })\overline{h} \big) {\pi }_{\ast }(X({\pi }^{\ast }\overline{f})) 
+ \big( ({\pi }_{\ast } \comp {\pi }^{\ast })\overline{f} \big) {\pi }_{\ast }(X({\pi }^{\ast }\overline{h})) \\
& = \overline{h} Y(\overline{f}) + \overline{f} Y(\overline{h}). 
\end{align*}
Since $Y$ is a linear mapping of $C^{\infty}(M/G)$ into itself, it follows that it is a 
derivation of $C^{\infty}(M/G)$. 
\par We now show that the map $X \rightarrow {\pi }_{\ast } \comp X \comp {\pi }^{\ast }$ is a module homomorphism. For $X$, $X' \in (\mathrm{Der}\, C^{\infty}(M))^G$ and 
$\overline{f} \in C^{\infty}(M/G)$ we have 
\begin{align*}
\big( {\pi }_{\ast } \comp (X+X') \comp {\pi }^{\ast} \big) \overline{f} & = 
{\pi }_{\ast } \big( X({\pi }^{\ast }\overline{f}) + X'({\pi }^{\ast }\overline{f}) \big) \\
& \hspace{-.75in} = {\pi }_{\ast }(X({\pi }^{\ast }\overline{f})) + 
{\pi }_{\ast }(X'({\pi }^{\ast }\overline{f})) = ({\pi }_{\ast }\comp X \comp {\pi }^{\ast })(\overline{f}) + 
({\pi }_{\ast }\comp X' \comp {\pi }^{\ast })(\overline{f}) . 
\end{align*}
Hence the map $X \rightarrow {\pi }_{\ast } \comp X \comp {\pi }^{\ast }$ is linear. 
For every $h \in C^{\infty}(M)^G$ 
\begin{align*}
({\pi }_{\ast } \comp (h X) \comp {\pi }^{\ast }) \overline{f} & = 
{\pi }_{\ast }(h X({\pi }^{\ast }\overline{f})) = {\pi }_{\ast }(h)\, 
{\pi }_{\ast }(X({\pi }^{\ast }\overline{f})) \\
& = {\pi }_{\ast }(h) ( {\pi }_{\ast }\comp X \comp {\pi }^{\ast })\overline{f}. 
\end{align*}
Therefore the map given by Equation (\ref{eq-s2onestarnw}) is a module homomorphism. 
\hfill $\square $ \medskip 

The importance of the module homomorphism (\ref{eq-s2onestarnw}) stems from the following 
result. \medskip 

\noindent \textbf{Proposition 3.} \textit{Since $M$ is a smooth manifold, $(\mathrm{Der}\, C^{\infty}(M))^G 
= \mathfrak{X}(M)^G$. So $X \in \mathfrak{X}(M)^G$ implies that 
$Y = {\pi }_{\ast }\comp X \comp {\pi }^{\ast } \in \mathfrak{X}(M/G)$.} \medskip 

\noindent \textbf{Proof.} Because the orbit space $M/G$ is locally closed and subcartesian, every maximal integral curve of $X$ projects under the $G$ orbit map to a maximal integral curve of $Y$. It follows that $Y$ is a smooth vector field on $M/G$, see proposition 3.2.6 on page 34 of \cite{sniatycki}. \hfill $\square $ \medskip 

The following example shows that not every derivation on $C^{\infty}(M/G)$ is a vector field on 
$M/G$. \medskip 

\noindent \textbf{Example 2.} Consider the ${\Z }_2$ action on $\R $ generated by $\zeta : \R \rightarrow \R: x \mapsto -x$. The algebra $C^{\infty}(\R  )^{{\Z}_2}$ of smooth ${\Z }_2$ invariant functions is generated by the polynomial $\sigma (x) = x^2$. The orbit map of the ${\Z }_2$ action is $\pi : \R \rightarrow \R/{\Z }_2 \subseteq \R : x \mapsto \sigma (x) = \sigma $. The derivation $\frac{\partial }{\partial \sigma }$ of $C^{\infty}(\R  )^{{\Z}_2}$ is \emph{not} a smooth vector field on $\R /{\Z}_2$, because its maximal integral curve ${\gamma }_{{\sigma }_0}$ 
starting at ${\sigma }_0 \in \R/ {\Z }_2$, given by 
${\gamma }_{{\sigma }_0}(t) = t+{\sigma }_0$, is defined on $[-{\sigma }_0, \infty)$, which is \emph{not} an open interval that contains $0$. \hfill $\square $  \medskip

Fix $m \in M$. Then $G_m = \{ g \in G \setrule \, {\Phi }_g(m) = m \} $ is the isotropy group of the action $\Phi $ at $m$. It is a compact subgroup of $G$, see Duistermaat and Kolk \cite{duistermaat-kolk}. Let $H$ be a compact subgroup of $G$. The set 
\begin{equation}
M_H = \{ m \in M \setrule \, G_m = H \}
\label{eq-s2oneA}
\end{equation}
is a submanifold of $M$, which is not necessarily connected. Hence its 
connected component are submanifolds. Connected components of $M_H$ are $H$ invariant submanifolds of $M$, see Duistermaat and Kolk 
\cite{duistermaat-kolk}. The conjugacy class in $G$ of a 
closed subgroup $H$ is denoted by $(H) = \{ gHg^{-1} \in G \setrule \, g \in G \} $ and is called a \emph{type}. The set 
\begin{equation}
M_{(H)} = \{ m \in M \setrule , G_m \in (H) \} = G \cdot M_H
\label{eq-s2oneb}
\end{equation}
is called an \emph{orbit type} $(H)$. Moreover, the $G$ invariance of $(H)$ implies that each 
connected component of $M_{(H)}$ is $G$ invariant. The orbit type $M_{(H)}$ is \linebreak associated to the type $(H)$. Let $(H_1)$ and $(H_2)$ be two types. Define the partial order $\le $ by the condition: $(H_1) \le (H_2)$ if some $G$ conjugate of $H_2$ is a subgroup of $H_1$. Since the orbit space $\overline{M}$ is connected, there is a unique maximal orbit type $M_{(K)}$. \medskip 

\noindent \textbf{Proposition 4.} \textit{The maximal orbit type $M_{(K)}$ is open and dense in $M$.} \medskip 

\noindent \textbf{Proof.} See page 118 of Duistermaat and Kolk \cite{duistermaat-kolk}. \hfill $\square $ \medskip 

\noindent \textbf{Proposition 5.} \textit{The orbit type $M_{(H)}$ is a smooth invariant submanifold of every vector field $X \in \mathfrak{X}(M)^G$.} \medskip 

\noindent \textbf{Proof.} Let $m \in M_{(H)}$. Then $G_m = gHg^{-1}$ for some $g \in G$. Let ${\varphi }_t$ be the local $1$ parameter group of local diffeomorphisms of $M$ 
generated by $X \in \mathfrak{X}(M)^G$. Then  
\begin{align}
{\Phi}_{gHg^{-1}}\big( {\varphi }_t(m) \big) & = {\varphi }_t\big( {\Phi }_{gHg^{-1}}(m) \big)  
= {\varphi }_t(m),  \notag 
\end{align}
since $gHg^{-1} = G_m$. Thus $gHg^{-1} \subseteq G_{{\varphi }_t(m)}$. Conversely, suppose that $k \in G_{{\varphi }_t(m)}$. Then ${\varphi }_t(m) = {\Phi }_k\big( {\varphi }_t(m) \big) = 
{\varphi }_t\big( {\Phi }_k(m) \big)$. So 
\begin{displaymath}
m = {\varphi }_{-t}\big( {\varphi }_t(m) \big) = 
{\varphi }_{-t}\big( {\varphi }_t({\Phi }_k(m) \big) = {\Phi }_k(m),
\end{displaymath}
that is, $k \in G_m = gHg^{-1}$. Thus $G_{{\varphi }_t(m)} \subseteq gHg^{-1}$. Consequently, $G_{{\varphi }_t(m)} = gHg^{-1}$, that is, ${\varphi }_t(m) \in M_{(H)}$. 
\hfill $\square $ \medskip 

Let $Y$ be a smooth vector field on $\overline{M}$, which is $\pi $ related to the smooth $G$ invariant vector field $X$ on $M$, that is, $T_m\pi X(m) = Y\big( \pi (m) \big) $ for  $m \in M$. Then 
${\overline{\varphi }}_t(\overline{m}) = \pi \big({\varphi }_t(m) \big) $ for $(t, \overline{m}) \in 
\mathcal{D} = \{ (t, \overline{m}) \in \R \times \overline{M} \setrule \, (t, m ) \in D \, \, 
\& \, \, \pi (m ) = \overline{m} \}$.  The set $\mathcal{D}$ is well defined, because 
$(t, g \cdot m) \in D$ for every $g \in G$. So $\mathcal{D}$ is the domain of a local generator of $Y$, that is, $\overline{\varphi} : \mathcal{D} \subseteq \R \times \overline{M} \mapsto \overline{M}$ is a differentiable mapping such that 
$\frac{\dee {\overline{\varphi }}_t}{\dee t} (\overline{m}) = 
Y({\overline{\varphi }}_t(\overline{m}))$ for $(t,\overline{m}) \in \mathcal{D}$. \medskip 

Note that the orbit type $M_{(H)}$ need not be connected and its connected components may be of different dimensions. In the following we concentrate our attention on the properties of the connected components of $M_{(H)}$, which we denote by ${\mathrm{M}}_{(H)}$. \medskip 

\noindent \textbf{Proposition 6.} \textit{For every compact subgroup $H$ of $G$ the image of each connected component ${\mathrm{M}}_{(H)}$ of 
$M_{(H)}$ under the orbit mapping $\pi : M \rightarrow M/G$ is a smooth 
submanifold of the differential space $\big( M/G, C^{\infty}(M/G) \big)$.} \medskip 

\noindent \textbf{Proof.} See page 74 of \cite{sniatycki}. \hfill $\square $ \medskip 

\noindent \textbf{Proposition 7.} \textit{The connected component 
${\overline{\mathrm{M}}}_{(H)} = \pi ({\mathrm{M}}_{(H)})$ of the orbit type 
${\overline{M}}_{(H)}$ of $M/G$ is an invariant manifold of every smooth vector field $Y$ on $M/G$.} \medskip 

\vspace{-.15in}\noindent \textbf{Proof.} Let $\overline{\varphi }: \overline{U} \subseteq \R \times M/G \rightarrow M/G$ be the local flow of the vector field $Y$. For each point $y \in {\overline{\mathrm{M}}}_{(H)}$ there is an open neighborhood $\overline{V}$ of $y$ in 
$M/G$ such that for every $t \in [0, \varepsilon )$ the map ${\overline{\varphi }}_t$ is a local diffeomorphism onto its image. Hence ${\overline{\varphi }}_t(y) \in 
{\overline{\mathrm{\mathrm{M}}}}_{(H)}$ for every $t \in [0,\varepsilon )$. So  
${\overline{\mathrm{\mathrm{M}}}}_{(H)}$ is an invariant manifold of the 
vector field $Y$. \hfill $\square $ \medskip  

\noindent \textbf{Theorem 1.} \textit{Let $H$ be a compact subgroup of the Lie group $G$. Let $Y$ be a smooth vector field on $M/G$. Then on every 
connected component ${\mathrm{M}}_{(H)}$ of the orbit type $M_{(H)}$ there is a smooth $G$ invariant vector field $X$, which is 
${\pi }_{|{\mathrm{M}}_{(H)}}$ related to $Y_{|{\overline{\mathrm{M}}}_{(H)}}$.} \medskip 

We need the next few results to prove this theorem, which is the main result of this section. \medskip 

\noindent \textbf{Lemma 1.} \textit{The $G$ orbit $G\cdot m = \mathcal{O}$ through $m \in M$ is a submanifold of $M$.} \medskip 

\noindent \textbf{Proof.} Let $S_m$ be a slice to the $G$ action $\Phi $ at $m$. By Bochner's lemma, see \cite[p.306]{cushman-bates}, there is a local diffeomorphism 
$\psi : T_mM \rightarrow M$, which sends $0_m \in T_mM$ to $m \in M$ and 
intertwines the $H = G_m$ action 
\begin{displaymath}
{\Psi}_m: H \times T_mM \rightarrow T_mM: (h,v_m) \mapsto h \star v_m = 
T_m{\Phi }_h\, v_m
\end{displaymath} 
on $T_mM$ with the $H$ action ${\Phi }_H: H \times M \rightarrow M: (h,m) \mapsto {\Phi }_h(m)$. Since $T_mS_m$ is $H$ invariant, it follows that $\psi : T_mS_m \rightarrow S_m$ is a local diffeomorphism which sends $0_m$ to $m$. Let $L$ 
be a complement of $\mathfrak{h}$ in $\mathfrak{g}$, where $\mathfrak{h}$ is the Lie algbebra of $H$. The map 
\begin{displaymath}
\phi : L \times T_mS_m \rightarrow M: (\xi , v_m) \mapsto 
{\Phi }_{\exp \xi} \big( \psi (v_m) \big) , 
\end{displaymath}
which sends $(0_L,0_m)$ to $m$ is a local diffeomorphism that sends an open neighborhood of $(0_L,0_m)$ in $L \times \{ 0_m \}$ onto an open neighborhood of $m$ in $\mathcal{O}$. Thus $\mathcal{O}$ is a smooth submanifold of $M$ near $m$. For every $g \in G$, since ${\Phi }_g$ is a diffeomorphism of $M$, the map ${\Phi }_g \comp \phi $ is a local diffeomorphism of $(0_L,0_m)$ in $L \times \{ 0_m \}$ onto an open neighborhood of 
$g \cdot m$ in $\mathcal{O}$. Thus $\mathcal{O}$ is a submanifold of $M$. \medskip 

\noindent \textbf{Lemma 2.} \textit{For each connected component 
${\mathrm{M}}_{(H)}$ of the orbit type $M_{(H)}$ the map 
\begin{equation}
{\pi }_{|{\mathrm{M}}_{(H)}}: {\mathrm{M}}_{(H)} \rightarrow {\overline{\mathrm{M}}}_{(H)} 
= \pi ({\mathrm{M}}_{(H)}): m \mapsto \overline{m} 
\label{eq-two}
\end{equation}
is a smooth surjective submersion, whose typical fiber is an orbit of the $G$ action $\Phi $ restricted to $G \times {\mathrm{M}}_{(H)}$.} \medskip %

\noindent \textbf{Proof.} The orbit map $\pi : M \rightarrow \overline{M} = M/G$ is a surjective smooth map of the smooth manifold $M$ onto the differential space 
$(M/G, C^{\infty}(M/G))$. Hence its restriction 
${\pi }_{|{\mathrm{M}}_{(H)}}$ to the connected component ${\mathrm{M}}_{(H)}$ of the orbit type $M_{(H)}$ and the codomain to ${\overline{\mathrm{M}}}_{(H)} = \pi ({\mathrm{M}}_{(H)})$ is a smooth map of the smooth manifold 
${\mathrm{\mathrm{M}}}_{(H)}$ onto the differential space 
$({\overline{\mathrm{M}}}_{(H)}, C^{\infty}({\overline{\mathrm{M}}}_{(H)}))$. By Proposition 6, the differential space $({\overline{\mathrm{M}}}_{(H)}, 
C^{\infty}({\overline{\mathrm{M}}}_{(H)}))$ is a smooth manifold. 
Hence ${\pi }_{|{\mathrm{M}}_{(H)}}: {\mathrm{M}}_{(H)} \rightarrow 
{\overline{\mathrm{M}}}_{(H)}$ is a smooth map of the smooth manifold 
${\mathrm{M}}_{(H)}$ onto the smooth manifold ${\overline{\mathrm{M}}}_{(H)}$. At $\overline{m} \in {\overline{\mathrm{M}}}_{(H)}$, the fiber 
$({\pi }_{|{\mathrm{M}}_{(H)}})^{-1}(\overline{m})$ is the $G$ orbit in ${\mathrm{M}}_{(H)}$ through $m$, which is a smooth submanifold of ${\mathrm{M}}_{(H)}$. We have $T_m{\mathrm{M}}_{(H)} = \ker T_m {\pi }_{|{\mathrm{M}}_{(H)}} \oplus (\ker T_m {\pi }_{|{\mathrm{M}}_{(H)}})^{\perp}$, using the restriction of the $G$ invariant Riemannian metric on $M$ to 
$\ker T_m {\pi }_{|{\mathrm{M}}_{(H)}}$, 
see Palais \cite{palais}. Because the vector space 
$(\ker T_m {\pi }_{|{\mathrm{M}}_{(H)}})^{\perp}$ is isomorphic to 
$T_{\overline{m}}{\overline{\mathrm{M}}}_{(H)}$, having the same dimension, it follows that the map $T_m{\pi }_{|{\mathrm{M}}_{(H)}}$ is surjective. Consequently, the map 
${\pi }_{|{\mathrm{M}}_{(H)}}$ is a submersion. \hfill $\square $ \medskip %

Next we construct a connection on the fibration 
${\pi }_{|{\mathrm{M}}_{(H)}}: {\mathrm{M}}_{(H)} \rightarrow {\overline{\mathrm{M}}}_{(H)}$ and then review some geometric facts about geodesics. Because the $G$ action $\Phi $ (\ref{eq-s0zero}) on the smooth manifold $M$ is proper, it has a $G$ invariant Riemannian metric. Let 
$\mathbf{k}$ be the restriction of this metric to the smooth submanifold 
${\mathrm{M}}_{(H)}$. For each $m \in {\mathrm{M}}_{(H)}$ this yields the $G$ invariant decomposition 
\begin{equation}
T_m{\mathrm{M}}_{(H)} = {\mathrm{ver}}_m \oplus {\mathrm{hor}}_m, 
\label{eq-clonedagger}
\end{equation}
where ${\mathrm{ver}}_m = \ker T_m{\pi }_{|{\mathrm{M}}_{(H)}}$ and ${\mathrm{hor}}_m = 
({\mathrm{ver}}_m)^{\perp}$, using the metric $\mathbf{k}(m)$ on $T_m{\mathrm{M}}_{(H)}$. 
The distributions $\mathrm{ver}: {\mathrm{M}}_{(H)} \rightarrow T{\mathrm{M}}_{(H)}: m \mapsto {\mathrm{ver}}_m$ and $\mathrm{hor}:{\mathrm{M}}_{(H)} \rightarrow T{\mathrm{M}}_{(H)}: 
m \mapsto {\mathrm{hor}}_m$ are smooth. Moreover, for every $m \in {\mathrm{M}}_{(H)}$ the map 
\begin{equation}
(T_m{\pi }_{|{\mathrm{M}}_{(H)}})_{|{\mathrm{hor}}_m}: {\mathrm{hor}}_m \subseteq T_m{\mathrm{M}}_{(H)} 
\rightarrow T_{\overline{m}}{\overline{\mathrm{M}}}_{(H)}: v_m \mapsto 
{\overline{v}}_{\overline{m}}, 
\label{eq-cltwodagger}
\end{equation}
where $\overline{m} = \pi (m)$, is an isomorphism of vector spaces. Thus Equations (\ref{eq-clonedagger}) and (\ref{eq-cltwodagger}) define an Ehresmann connection $\mathcal{E}$ on the fibration 
${\pi }_{|{\mathrm{M}}_{(H)}}: {\mathrm{M}}_{(H)} \rightarrow {\overline{\mathrm{M}}}_{(H)}$. Because 
\begin{displaymath}
{\mathrm{ver}}_{g \cdot m} \oplus {\mathrm{hor}}_{g \cdot m}  = 
T_m{\Phi }_g (T_m M^{\, c}_{(H)}) = T_m{\Phi }_g {\mathrm{ver}}_m \oplus 
T_m{\Phi }_g{\mathrm{hor}}_m 
\end{displaymath}
and $T_m{\Phi }_g{\mathrm{ver}}_m \subseteq {\mathrm{ver}}_{g\cdot m}$ imply 
${\mathrm{ver}}_{g \cdot m} = T_m{\Phi }_g {\mathrm{ver}}_m$ and 
${\mathrm{hor}}_{g\cdot m} = T{\Phi }_g {\mathrm{hor}}_m$, the distributions 
$\mathrm{ver}$ and $\mathrm{hor}$ are $G$ invariant. Thus the connection 
$\mathcal{E}$ is $G$ invariant. \medskip 

Let $\tau : T{\mathrm{M}}_{(H)} \rightarrow {\mathrm{M}}_{(H)}$ and 
$\rho : T^{\ast}{\mathrm{M}}_{(H)} \rightarrow {\mathrm{M}}_{(H)}$ be the 
tangent and cotangent bundle projection maps, respectively. The metric 
$\mathbf{k}$ on ${\mathrm{M}}_{(H)}$ defines a vector bundle isomorphism 
\begin{displaymath}
{\mathbf{k}}^{\sharp}: T{\mathrm{M}}_{(H)} \rightarrow T^{\ast }{\mathrm{M}}_{(H)}: v_m \mapsto p_m = 
{\mathbf{k}}^{\sharp}(m)(v_m), 
\end{displaymath}
where $\langle {\mathbf{k}}^{\sharp}(m)(v_m) | v_m \rangle = \mathbf{k}(m)(v_m,v_m)$.
The inverse of ${\mathbf{k}}^{\sharp}$ is ${\mathbf{k}}^{\flat }$. 
The metric $\mathbf{k}$ determines the Hamiltonian function 
\begin{equation}
E: T^{\ast }{\mathrm{M}}_{(H)} \rightarrow \R : p_m \mapsto 
\onehalf \mathbf{k}(m)\big( {\mathbf{k}}^{\flat }(m)(p_m), {\mathbf{k}}^{\flat }(m)(p_m) \big) ,
\label{eq-clthree}
\end{equation}
which gives rise to the Hamiltonian system $(E, T^{\ast }{\mathrm{M}}_{(H)}, \omega )$, 
where $\omega $ is the canonical symplectic form on 
$T^{\ast }{\mathrm{M}}_{(H)}$. The 
Hamiltonian vector field $X_E$ on $T^{\ast }{\mathrm{M}}_{(H)}$ is defined by 
$X_E \lefthook \, \omega = \dee E$. For $p_m \in T^{\ast }_m{\mathrm{M}}_{(H)}$ let 
\begin{displaymath}
{\varphi }^{X_E}_t: T^{\ast }{\mathrm{M}}_{(H)} \rightarrow 
T^{\ast }{\mathrm{M}}_{(H)}: p_m \mapsto 
{\varphi }^{X_E}_t(p_m)
\end{displaymath}
be the local flow of the vector field $X_E$, which is defined for $t$ in an open interval $I_m$ in $\R $ containing $0$. For $v_m \in T_m{\mathrm{M}}_{(H)}$ the curve ${\gamma }_{v_m}$ given by $t \rightarrow \big( \rho \comp {\varphi }^{X_E}_t \big) ({\mathbf{k}}^{\sharp}(m)(v_m))$ 
is a \emph{geodesic} on ${\mathrm{M}}_{(H)}$, starting at $m \in 
{\mathrm{M}}_{(H)}$, for the metric $\mathbf{k}$. There is an open tubular neighborhood $U$ of the zero section of the cotangent bundle $\rho :T^{\ast }{\mathrm{M}}_{(H)} \rightarrow {\mathrm{M}}_{(H)}$ 
such that the local flow ${\varphi }^{X_E}_t$ is defined for all $t \in [0,1]$. 
For each $m \in {\mathrm{M}}_{(H)}$, the \emph{exponential map} 
\begin{equation}
\begin{array}{l}
{\mathrm{exp}}_m: U_{0_m} \subseteq T_m{\mathrm{M}}_{(H)} \rightarrow {\mathrm{M}}_{(H)}: \\
\rule{0pt}{12pt}\hspace{.45in}v_m \mapsto {\gamma }_{v_m}(1) = 
\big( \rho \comp {\varphi }^{X_{E}}_1\big) ({\mathbf{k}}^{\sharp}(m)(v_m)) ,
\end{array} 
\label{eq-clone}
\end{equation}
is a diffeomorphism onto $V_m = {\mathrm{exp}}_m(U_{0_m})$, where $U_{0_m} \subseteq U$ is a suitable open neighborhood of $0_m \in T_m{\mathrm{M}}_{(H)}$, see Brickell and Clark \cite{brickell-clark}. \medskip 

Next we reduce the $G$ symmetry of the Hamiltonian system 
$(E, T^{\ast }{\mathrm{M}}_{(H)}, \omega )$. 
Because the metric $\mathbf{k}$ on ${\mathrm{M}}_{(H)}$ is $G$ invariant, the smooth Hamiltonian $E$ (\ref{eq-clthree}) on ${\mathrm{M}}_{(H)}$ is $G$ invariant, it induces a metric 
$\overline{\mathbf{k}}$ on ${\overline{\mathrm{M}}}_{(H)}$ such that 
${\pi }^{\ast }\overline{\mathbf{k}} = \mathbf{k} \comp (\mathrm{perp}, \mathrm{perp})$, where 
$\mathrm{perp}: T{\mathrm{M}}_{(H)} \rightarrow (T{\mathrm{M}}_{(H)})^{\perp}$ is orthogonal projection. The smooth Hamiltonian $E$ (\ref{eq-clthree}) on ${\mathrm{M}}_{(H)}$ is $G$ invariant, and hence induces a smooth Hamiltonian function 
\begin{displaymath}
\overline{E}: T^{\ast }{\overline{\mathrm{M}}}_{(H)} \rightarrow \R :{\overline{p}}_{\overline{m}} 
\mapsto \onehalf \overline{\mathbf{k}}(\overline{m}) 
\big( {\overline{\mathbf{k}}}^{\, \flat}(\overline{m})({\overline{p}}_{\overline{m}} ), 
{\overline{\mathbf{k}}}^{\, \flat}(m)({\overline{p}}_{\overline{m}} ) \big) . 
\end{displaymath}
Since the $G$ orbit map ${\pi }_{|{\mathrm{M}}_{(H)}}: {\mathrm{M}}_{(H)} \rightarrow 
{\overline{\mathrm{M}}}_{(H)}$ is smooth, symplectic reduction of the Hamiltonian 
system $(E, T^{\ast }{\mathrm{M}}_{(H)}, \omega )$ leads to the $G$ reduced Hamiltonian 
system $(\overline{E}, T^{\ast }{\overline{\mathrm{M}}}_{(H)}, \overline{\omega })$. The reduced system has a Hamiltonian vector field $X_{\overline{E}}$ defined by $X_{\overline{E}} \lefthook \overline{\omega } = \dee \overline{E}$. Its local flow 
\begin{displaymath}
{\varphi }^{X_{\overline{E}}}_t: T^{\ast }{\overline{\mathrm{M}}}_{(H)} \rightarrow 
T^{\ast }{\overline{\mathrm{M}}}_{(H)}: {\overline{p}}_{\overline{m}} \mapsto 
{\varphi }^{X_{\overline{E}}}_t({\overline{p}}_{\overline{m}}) 
\end{displaymath}
is $\pi $ related to the local flow ${\varphi }^{X_E}_t$ of $X_E$, that is, 
$\pi \comp {\varphi }^{X_E}_t = {\varphi }^{X_{\overline{E}}}_t \comp \pi $. The curve ${\overline{\gamma }}_{{\overline{v}}_{\overline{m}}}$ given by 
$t \mapsto \big( \overline{\rho } \comp {\varphi }^{X_{\overline{E}}}_t\big) 
({\overline{\mathbf{k}}}^{\sharp}(\overline{m})({\overline{v}}_{\overline{m}}) \big) $ is a geodesic on ${\overline{\mathrm{M}}}_{(H)}$ starting at $\overline{m}$ for 
the reduced metric $\overline{\mathbf{k}}$. Here 
$\overline{\rho }:T^{\ast }{\overline{\mathrm{M}}}_{(H)} 
\rightarrow {\overline{\mathrm{M}}}_{(H)}$ is the cotangent bundle projection map. Note that $\pi \comp {\gamma }_{v_m} = 
{\gamma }_{{\overline{v}}_{\overline{m}}} \comp \pi $, 
where ${\overline{v}}_{\overline{m}} = T_m\pi \, v_m$. There is an open neighborhood $\overline{U} = T\pi (U)$ of the zero section of $\overline{\rho }:
T^{\ast }{\overline{\mathrm{M}}}_{(H)} 
\rightarrow {\overline{\mathrm{M}}}_{(H)}$ such that the local flow ${\varphi }^{X_{\overline{E}}}_t$ of the reduced vector field $X_{\overline{E}}$ is defined for all $t \in [0,1]$. The reduced exponential map
\begin{displaymath}
\begin{array}{l}
{\overline{\mathrm{exp}}}_{\, \overline{m}}: 
{\overline{U}}_{{\overline{0}}_{\overline{m}}} 
\subseteq T_{\overline{m}}{\overline{\mathrm{\mathrm{M}}}}_{(H)} \rightarrow 
{\overline{\mathrm{M}}}_{(H)} : \\
\hspace{.55in} {\overline{v}}_{\overline{m}} \mapsto 
{\gamma }_{{\overline{v}}_{\overline{m}}}(1) = 
\big( \overline{\rho }\comp {\varphi }^{X_{\overline{E}}}_t \big) 
({\overline{\mathbf{k}}}^{\flat}({\overline{v}}_{\overline{m}}) )
\end{array}
\end{displaymath}
is a diffeomorphism onto ${\overline{V}}_{\overline{m}} = 
{\overline{\mathrm{exp}}}_{\overline{m}}({\overline{U}}_{{\overline{0}}_{\overline{m}}})$, where 
${\overline{U}}_{{\overline{0}}_{\overline{m}}} = \pi (U_{0_m})$. \medskip 

\noindent \textbf{Proposition 8.} \textit{The fibration ${\pi }_{|{\mathrm{M}}_{(H)}}: {\mathrm{M}}_{(H)} \rightarrow {\overline{\mathrm{M}}}_{(H)} = 
\pi ({\mathrm{M}}_{(H)})$ is locally trivial.} \medskip 

\noindent \textbf{Proof.} For some $b_m > 0$, the open 
ball $B_m = \{ v_m \in {\mathrm{hor}}_m \setrule \, \mathbf{k}(m)(v_m,v_m) < b_m \}$ is a subset of $U_{0_m}$. Then $V_m = {\mathrm{exp}}_m(B_m)$ is a submanifold of ${\mathrm{M}}_{(H)}$ containing $m$. Look at the geodesic ${\gamma }_{T_m{\Phi }_gv_m} 
= {\gamma }_{v_{g\cdot m}}$ on ${\mathrm{M}}_{(H)}$ given by 
\begin{displaymath}
t \mapsto \big( \rho \comp {\varphi }^{X_E}_t \big) 
({\mathbf{k}}^{\sharp}(g \cdot m)(T_m{\Phi }_g v_m))
\end{displaymath} 
starting at $g \cdot m$. One has $B_{g \cdot m} = T_m{\Phi }_g(B_m)$. To see this, observe that there is a $v_m \in T_m {\mathrm{M}}_{(H)}$ such that $v_{g\cdot m} = 
T_m{\Phi }_g v_m$, since ${\Phi }_g$ is a diffeomorphism. So 
\begin{align*}
b_{g\cdot m} & = \mathbf{k}(g\cdot m)(v_{g\cdot m}, v_{g \cdot m}) = 
\mathbf{k}(g\cdot m)(T_m{\Phi }_gv_m, T_m{\Phi }_gv_m) \\
& = \mathbf{k}(m)(v_m,v_m) = b_m. 
\end{align*}
Thus ${\mathrm{exp}}_{g\cdot m} = {\Phi }_g \comp {\mathrm{exp}}_m$ is a diffeomorphism 
of the open ball $B_{g\cdot m}$ of \linebreak
radius $b_{g\cdot m} = b_m$ contained in $U_{g\cdot m}$ 
onto a submanifold $V_{g\cdot m} = {\mathrm{exp}}_{g\cdot m}(B_{g\cdot m}) = 
{\Phi }_g(V_m)$. \medskip 

For every ${\overline{m}}' \in {\overline{V}}_{\overline{m}}$ let 
\begin{displaymath}
{\overline{\gamma }}_{\overline{m}, {\overline{m}}'}: [0,1] \rightarrow {\overline{V}}_{\overline{m}}: 
t \mapsto {\overline{\exp }}_{\overline{m}} \, t{\overline{v}}'_{\overline{m}}
\end{displaymath}
be the geodesic in ${\overline{V}}_{\overline{m}}$ joining $\overline{m}$ to 
${\overline{m}}'$, that is, ${\overline{\gamma }}_{\overline{m}, {\overline{m}}'}(1) = 
{\overline{m}}'$. Because ${\overline{\exp }}_{\overline{m}}$ is a diffeomorphism, 
the vector ${\overline{v}}'_{\overline{m}} \in B_{\overline{m}} = \pi (B_m)$ is uniquely 
determined by ${\overline{m}}'$. Let 
\begin{displaymath}
{\mathcal{P}}_{{\overline{\gamma }}_{\overline{m}, {\overline{m}}'}}: 
({\pi }_{|{{\mathrm{M}}_{(H)}}})^{-1}(\overline{m}) \rightarrow 
({\pi }_{|{{\mathrm{M}}_{(H)}}})^{-1}({\overline{m}}') : n \rightarrow {\gamma }_{n,n'}(1) = n', 
\end{displaymath}
where ${\gamma }_{n,n'}: [0,1] \rightarrow V_n: t \mapsto {\exp }_n tv'_n$ 
is the horizontal lift of the geodesic ${\overline{\gamma }}_{ \overline{m}, {\overline{m}}' }$ using 
the connection $\mathcal{E}$. Here 
$n' \in B_n$ with $\pi (n') = {\overline{m}}'$, $n = g \cdot m$, and $v'_n \in B_n$. The map 
${\mathcal{P}}_{{\overline{\gamma }}_{\overline{m}, {\overline{m}}'}}$ is 
\emph{parallel translation} of the fiber $({\pi }_{|{{\mathrm{M}}_{(H)}}})^{-1}(\overline{m})$ 
along the geodesic ${\overline{\gamma }}_{\overline{m}, {\overline{m}}'}$ 
joining $\overline{m}$ to ${\overline{m}}'$ in ${\overline{V}}_{\overline{m}}$ using the connection $\mathcal{E}$. \medskip 

Consider the mappings 
\begin{displaymath}
{\tau }_{{\overline{V}}_{\overline{m}} }: {\overline{V}}_{\overline{m}} \times 
({\pi }_{|{{\mathrm{M}}_{(H)}}})^{-1}(\overline{m}) \rightarrow 
({\pi }_{|{{\mathrm{M}}_{(H)}}})^{-1}({\overline{V}}_{\overline{m}}): 
({\overline{m}}', n) \mapsto {\mathcal{P}}_{{\overline{\gamma }}_{ \overline{m}, {\overline{m}}' }}(n) 
\end{displaymath}
and the projection mapping 
\begin{displaymath}
{\pi }_1 = ({\pi }_{ |{\mathrm{M}}_{(H)}})_
{|({\pi }_{|{\mathrm{M}}_{(H)}})^{-1}({\overline{V}}_{\overline{m}}) }: 
({\pi }_{|{\mathrm{M}}_{(H)}})^{-1}({\overline{V}}_{\overline{m}}) \rightarrow 
{\overline{V}}_{\overline{m}}: n \mapsto {\pi }_{| {\mathrm{M}}_{(H)} }(n). 
\end{displaymath}
Then ${\tau }_{{\overline{V}}_{\overline{m}}}$ is a local trivialization of the fibration 
defined by the mapping ${\pi }_1$ because for every $n \in 
({\pi }_{|{{\mathrm{M}}_{(H)}}})^{-1}(\overline{m})$ and every ${\overline{m}}' \in 
{\overline{V}}_{\overline{m}}$ 
\begin{displaymath}
({\pi }_1 \comp {\tau }_{{\overline{V}}_{\overline{m}}})({\overline{m}}', n) = 
{\pi }_1\big( {\mathcal{P}}_{{\overline{\gamma }}_{\overline{m}, {\overline{m}}'}}(n) \big) = {\overline{m}}'. 
\end{displaymath} 
We now show that ${\tau }_{{\overline{V}}_{\overline{m}}}$ is a diffeomorphism. Define the 
smooth maps
\begin{displaymath}
\rho : {\overline{V}}_{\overline{m}} \times 
({\pi }_{|{\mathrm{M}}_{(H)}})^{-1}({\overline{m}}) \rightarrow 
({\pi }_{|{\mathrm{M}}_{(H)}})^{-1}({\overline{V}}_{\overline{m}}): 
({\overline{m}}', n) \mapsto 
{\mathcal{P}}_{-{\overline{\gamma }}_{\overline{m}, {\overline{m}}'}}(n)
\end{displaymath}
and 
\begin{displaymath}
\sigma : ({\pi }_{|{\mathrm{M}}_{(H)}})^{-1}({\overline{V}}_{\overline{m}}) \rightarrow 
{\overline{V}}_{\overline{m}} \times ({\pi }_{|{\mathrm{M}}_{(H)}})^{-1}({\overline{m}}): 
n \mapsto \big( {\pi }_1(n), \rho ({\pi }_1(n), n) \big) .
\end{displaymath}
The following calculation shows that $\sigma \comp {\tau }_{{\overline{V}}_{\overline{m}}} = 
{\mathrm{id}}_{{\overline{V}}_{\overline{m}} \times 
({\pi }_{|{\mathrm{M}}_{(H)}})^{-1}({\overline{m}})}$. 
\begin{align*}
(\sigma \comp {\tau }_{{\overline{V}}_{\overline{m}}})({\overline{m}}', n) & = 
\sigma \big( {\mathcal{P}}_{{\overline{\gamma }}_{ \overline{m}, {\overline{m}}' }}(n) \big) \\
& = \big( {\pi }_1({\mathcal{P}}_{{\overline{\gamma }}_{ \overline{m}, {\overline{m}}' }}(n)), 
\rho ({\pi }_1({\mathcal{P}}_{{\overline{\gamma }}_{ \overline{m}, {\overline{m}}' }}(n)), 
{\mathcal{P}}_{{\overline{\gamma }}_{\overline{m}, {\overline{m}}' }}(n) ) \big) \\
& = \big( {\overline{m}}', \rho ({\overline{m}}', 
{\mathcal{P}}_{{\overline{\gamma }}_{ \overline{m}, {\overline{m}}' }}(n)) \big) 
= \big( {\overline{m}}', ({\mathcal{P}}_{-{\overline{\gamma }}_{ \overline{m}, {\overline{m}}' }} \comp 
{\mathcal{P}}_{{\overline{\gamma }}_{\overline{m}, {\overline{m}}' }})(n) \big) \\
& = ({\overline{m}}', n) .  
\end{align*}
Also 
\begin{align*}
({\tau }_{{\overline{V}}_{\overline{m}}} \comp \sigma )(n) & = 
{\tau }_{{\overline{V}}_{\overline{m}}}\big( {\pi }_1(n), \rho ({\pi }_1(n),n) \big) 
= ({\mathcal{P}}_{{\overline{\gamma }}_{ \overline{m}, {\overline{m}}' }} \comp 
{\mathcal{P}}_{ -{\overline{\gamma }}_{\overline{m}, {\overline{m}}' }})(n) = n, 
\end{align*}
that is, ${\tau }_{{\overline{V}}_{\overline{m}}} \comp \sigma = 
{\mathrm{id}}_{{\pi }_{|{\mathrm{M}}_{(H)}})^{-1}({\overline{V}}_{\overline{m}}) }$. 
Thus ${\tau }_{{\overline{V}}_{\overline{m}}}$ is a diffeomorphism. \medskip 
 
The preceding argument can be repeated at each point of ${\mathrm{M}}_{(H)}$. Hence the fibration ${\pi }_{|{\mathrm{M}}_{(H)}}: {\mathrm{M}}_{(H)} \rightarrow 
{\overline{\mathrm{M}}}_{(H)}$ is locally trivial. 
\hfill $\square $ \medskip 

\noindent \textbf{Corollary 1.} \textit{The locally trivial fibration defined by 
$\pi _{|{\mathrm{M}}_{(H)}}$ (\ref{eq-two}) has a local trivialization 
\begin{equation}
\begin{array}{l}
{\tau }_V: V \subseteq G \cdot {\mathrm{M}}_H = {\mathrm{M}}_{(H)} \rightarrow (U = {\pi }_{|{\mathrm{M}}_H)}(V) ) \times G: \\
\hspace{.5in} {\Phi }_k(g' \cdot m') \longmapsto ({\overline{m}}', kg') ,
\end{array} 
\label{eq-three}
\end{equation}
where $V$ is an open $G$ invariant subset of ${\mathrm{M}}_{(H)}$ with $k$, $g' \in G$ and $m' \in {\mathrm{M}}_H$.} \medskip 

\noindent \textbf{Proof.}  Suppose that $m'' = g' \cdot m'$ for some $m'$, $m'' \in 
{\mathrm{M}}_H$ and some $g' \in G$. 
Then 
\begin{align*}
({\pi }_{|{\mathrm{M}}_{(H)}}(m'), e) & = ({\pi }_{|{\mathrm{M}}_{(H)}}(m''), e) = 
{\tau }_V(e \cdot m'') = {\tau }_V(m'')   \\
&= {\tau }_V(g' \cdot m') = ({\pi }_{|{\mathrm{M}}_{(H)}}(m'), g'). 
\end{align*}
So $g' =e$. In other words, the $G$ orbit of $m'$ in ${\mathrm{M}}_H$ is 
$\{ m ' \}$. Hence ${\pi }_{|{\mathrm{M}}_{(H)}}(m') = \{ m' \} $ for every $m' \in 
{\mathrm{M}}_H$. Thus the map ${\tau }_V$ (\ref{eq-three}) is given by 
\begin{equation}
{\tau }_V: V \subseteq G \cdot {\mathrm{M}}_H \rightarrow U \times G: {\Phi }_k(g' \cdot m') \mapsto (\{ m' \} , kg') . 
\label{eq-five}
\end{equation}
Clearly, ${\tau }_V$ is a diffeomorphism. It intertwines the $G$ action $\Phi $ with the $G$ action 
\begin{equation}
\phi : G \times (U \times G) \rightarrow G \times U: \big(  k, (u,g) \big) \mapsto (u, kg) 
\label{eq-six}
\end{equation}
and satisfies ${\pi }_{|V} = {\pi }_1 \comp {\tau }_V$, where 
${\pi }_1: U \times G \rightarrow U: (u,g) \mapsto u$. 
Hence ${\tau }_V$ is a local trivialization. \hfill $\square $ \medskip 

\noindent \textbf{Lemma 3.} \textit{Every smooth vector field $Y_U$ on $U= 
{\pi }_{|{\mathrm{M}}_{(H)}}(V)$ is ${\pi }_1$ related to a smooth $G$ invariant vector field $\widetilde{Y}$ on $U \times G$.} \medskip 

\noindent \textbf{Proof.} To see this let $\xi \in \mathfrak{g}$, the Lie algebra of $G$. Let $\widetilde{Y}(u,g) = \big( Y_U(u), T_eL_g\xi \big)$, where $L_g$ is left translation on $G$ by $g$. By construction $\widetilde{Y}(u,g) \in T_uU \times T_gG = T_{(u,g)}(U \times G)$. So $\widetilde{Y}$ is a vector field on $U \times G$, which is smooth. For every $h \in G$ and every 
$(u,g) \in U\times G$ one has 
\begin{align*}
\widetilde{Y}\big( {\phi }_h(u,g) \big) & = \widetilde{Y}(u, hg) = \big( Y_U (u), T_eL_{hg} \xi \big) \\
& = T_{(u,g)} {\phi }_h ( Y_U(u), T_eL_g \xi ) = T_{(u,g)}{\phi }_h\, \widetilde{Y}(u,g).
\end{align*}
So $\widetilde{Y}$ is a $G$ invariant vector field on $U \times G$. Moreover, 
\begin{displaymath}
T_{(u,g)}{\pi }_1 \widetilde{Y}(u,g) = Y_U(u) = Y_U({\pi }_1(u,g)), \, \,  \mbox{for every $(u,g) \in U \times G$.}
\end{displaymath}
So the vector field $\widetilde{Y}$ on $U \times G$ and the vector field $Y_U$ on $U$ are ${\pi }_1$ related. \hfill $\square $ \medskip 

\noindent \textbf{Lemma 4.} \textit{Every smooth vector field $Y_U$ on $U$ is 
${\pi }_{|V}$ related to a smooth $G$ invariant vector field $X$ on $V$.} \medskip 

\noindent \textbf{Proof.} Pull the vector field $\widetilde{Y}$ on $U \times G$ back by the trivialization ${\tau }_V$ (\ref{eq-three}) to a vector field $X$ on $V$. Since ${\tau }_V$ intertwines the $G$ action $\Phi $ on $V$ with the $G$ action $\phi $ on $U \times G$, the vector field $X$ is $G$ invariant. For $m' \in V \cap {\mathrm{M}}_H$ and $g \cdot m' \in V$ one has
\begin{align*}
T_{g \cdot m' } {\pi }_{|V} \, X (g \cdot m' ) & = T_{g \cdot  m'} {\pi }_{|V} \big( T_{ {\tau }_V(g\cdot m')} {\tau }^{-1}_V 
\widetilde{Y}({\tau }_V(g\cdot m')) \big) \\
& = T_{g \cdot m'} \big( {\pi }_{|V} \comp {\tau }^{-1}_V\big) 
\widetilde{Y}( {\pi }_{|{\mathrm{M}}_H}(m') , g) \\
&= T_{({\pi }_{|{\mathrm{M}}_H}(m'),g)}{\pi }_1 \, \widetilde{Y}({\pi }_{|{\mathrm{M}}_H}(m'), g) = 
Y_U({\pi }_{|V}(g \cdot m')) . 
\end{align*}
Thus the $G$ invariant vector field $X$ on $V$ is ${\pi }_{|V}$ related to 
the vector field $Y_U$ on $U$. \hfill $\square $ \medskip 

\noindent \textbf{Proof of theorem 1.} We just have to piece the local bits together. Cover the orbit type ${\mathrm{M}}_{(H)}$ by 
${\{ (V_i, {\tau }_{V_i}) \} }_{i\in I}$, where 
\begin{displaymath}
{\tau }_{V_i}: V_i \subseteq G\cdot {\mathrm{M}}_H \rightarrow U_i \times G: g' \cdot m' \mapsto 
\big( {\pi }_{|V_i}(m'), g' \big) 
\end{displaymath}
is a local trivialization of the bundle ${\pi }_{|{\mathrm{M}}_{(H)}}: {\mathrm{M}}_{(H)} \rightarrow 
{\overline{\mathrm{M}}}_{(H)}$. Let $Y$ be a smooth vector field on 
${\overline{\mathrm{M}}}_{(H)} \subseteq \overline{M}$. Since 
${\pi }_{|{\mathrm{M}}_{(H)}}( V_i \cap {\mathrm{M}}_H) =U_i$ and 
${\pi }_{|{\mathrm{M}}_{(H)}}$ is an open mapping, $U_i$ is an open subset of 
${\overline{\mathrm{M}}}_{(H)}$. Because ${\{ V_i \}}_{i \in I}$ covers ${\mathrm{M}}_{(H)}$, it follows that ${\{ U_i \}}_{i \in I}$ is an open covering of ${\overline{\mathrm{M}}}_{(H)}$. Applying Lemma 3 to the smooth vector field $Y_{U_i} = Y|U_i$ and then using Lemma 4, we obtain a $G$ invariant vector field $X_{V_i}$ on $V_i$, which is ${\pi }_{|V_i}$ related to the vector field $Y_{U_i}$ on $U_i$. Since $Y$ is a smooth vector field on ${\overline{\mathrm{M}}}_{(H)}$, on $U_i \cap U_j$, where $i$, $j \in I$, one has $Y_{U_i} = Y_{U_j}$. So on $V_i \cap V_j$ one has $X_{V_i} = X_{V_j}$. Thus the $G$ invariant vector fields $X_{V_i}$ piece together to give a smooth $G$ invariant vector field $X$ on ${\mathrm{M}}_{(H)}$. Since $X_{V_i}$ is ${\pi }_{|V_i}$ related to the vector field $Y_{U_i}$, the vector field $X$ on 
${\mathrm{M}}_{(H)}$ is ${\pi }_{|{\mathrm{M}}_{(H)}}$ related to the vector field $Y$ on ${\overline{\mathrm{M}}}_{(H)}$. \hfill $\square $

\section{Vector fields on $\mathbf{M/G}$}

We start with a local argument in a neighbourhood of a point $m\in M$ with
compact isotropy group $H$. By Bochner's lemma, there is a local
diffeomorphism $\psi :T_{m}M\rightarrow M,~$which sends $0_{m}\in T_{m}M$ to 
$m\in M$ and intertwines the linear $H$ action 
\begin{equation}
\Psi :H\times T_{m}M\rightarrow T_{m}M:(h,v_{m})\mapsto h\cdot
v_{m}=T_{m}\Phi _{h}v_{m}  
\label{16}
\end{equation}%
with the $H$ action $\Phi _{H}:H\times M\rightarrow M:(h,m)\mapsto 
{\Phi }_{h}(m)$.\medskip

Because the action $\Phi $ on the smooth manifold $M$ is proper, it has a 
$G$ invariant Riemannian metric $\mathbf{k}.$ Using the restriction of 
$\mathbf{k}$ to $T_{m}M$, we define $E=T_{m}(G\cdot m)^{\perp }\subseteq T_{m}M$. Then, there is an $H$ invariant open ball $B\subseteq E$, centered at $0_{m}$ with $B$ contained in the domain of local diffeomorphism $\psi $ such that $\psi (B)=S_{m}$ is a slice to the $G$ action on $M$ at $m.$\medskip

We now construct \ a model for the $H$ orbit space $B/H$ of the restriction
to $B$ of the linear action $\Psi _{m}$ of $H$ on $T_{m}M.$ Let $%
\{v_{i}\}_{i=1}^{n}$ be a basis of the vector space $E$. Hence, $E$ is
isomorphic to $\mathbb{R}^{n}$. Let $\mathbf{x}=(x_{1},\ldots ,x_{n})$
coordinates on $\mathbb{R}^{n}$ with respect to the basis $\{v_{i}\}_{i=1}^{n}.$ Let $\{\sigma _{j}\}_{j=1}^{\ell }$ be a set of generators
for the algebra of $H$ invariant polynomials on $\mathbb{R}^{n}.$ Let 
$\mathbf{y}=(y_{1},...,y_{\ell })$ be coordinates on $\mathbb{R}^{\ell}$. The orbit
map of of the $H$ action on $\mathbb{R}^{n}$ is%
\begin{equation}
\rho :\mathbb{R}^{n}\rightarrow \mathbb{R}^{n}/H:\mathbf{x}\mapsto \mathbf{y}%
=(\sigma _{1}(\mathbf{x),\ldots ,\sigma }_{\ell }(\mathbf{x)).}  
\label{17}
\end{equation}%

By Schwarz's theorem $f\in C^{\infty }(\mathbb{R}^{n})$ is $H$ invariant,
i.e. $f\in C^{\infty }(\mathbb{R}^{n})^{H}$, if and only if there exists a
function $F\in C^{\infty }(\mathbb{R}^{l})$ such that $f(\mathbf{x})= 
F(\sigma _{1}(\mathbf{x),\ldots ,\sigma }_{\ell }(\mathbf{x))}$ for every 
$\mathbf{x}\in \mathbb{R}^{n}.$ Smooth functions on the orbit space 
$\mathbb{R}^{n}/H$ are restrictions to $\mathbb{R}^{n}/H$ of smooth functions on $\mathbb{R}^{l}$. For every $\overline{f}\in C^{\infty }(\mathbb{R}^{n}/H),$ the
pull-back $\rho ^{\ast }\overline{f}$ by the orbit map is given by
\[
\rho ^{\ast }\overline{f}(\mathbf{x})=(\overline{f}\comp \rho )(\mathbf{x})=
\overline{f}(\sigma _{1}(\mathbf{x}), \ldots ,\sigma _{\ell }(\mathbf{x})). 
\]

The $H$ orbit map $\rho :\mathbb{R}^{n}\rightarrow \mathbb{R}^{n}/H\subseteq 
\mathbb{R}^{\ell }$ is a smooth map of differential spaces. Consider an open
ball $B\subseteq \mathbb{R}^{n},$ centered at the origin $\mathbf{0}\in $ 
$\mathbb{R}^{n},$ with closure $\overline{B}.$ We are interested in $B/H,$ the
space of $H$ orbits in $B$. Restricting $\rho $ to the domain 
$B\subseteq \mathbb{R}^{n}$ and the codomain $\Sigma =\rho (B)\subseteq 
\mathbb{R}^{\ell }$ gives 
\begin{equation}
\rho _{B}:B\rightarrow \Sigma :\mathbf{x}\mapsto (\sigma _{1}(\mathbf{x),
\ldots ,\sigma }_{\ell }(\mathbf{x)),}  
\label{18}
\end{equation}%
which is a surjective smooth map of differential spaces. For technical
reasons, we also need to keep in mind an intermediate
restriction of the orbit map $\rho $ to the closure $\overline{B}$ of $B$, 
\[
\rho _{\, \overline{B}}:\overline{B}\rightarrow \overline{\Sigma }:\mathbf{x}\mapsto
(\sigma _{1}(\mathbf{x),\ldots ,\sigma }_{\ell }(\mathbf{x)),}
\]%
where $\overline{\Sigma}=\rho (\overline{B}).$ The following diagram, in which the horizontal arrows denote the inclusion map, 
\[
\begin{array}{ccccc}
&  &  &  &  \\ 
& B & \hookrightarrow  & \overline{B} &  \\ 
\rho _{B} & \downarrow  &  & \downarrow  & \rho _{\, \overline{B}} \\ 
& \Sigma  & \hookrightarrow  & \overline{\Sigma} &  \\ 
&  &  &  & 
\end{array}%
\]%
commutes and is invariant under the action of $H$ on $\mathbb{R}^{n}.$%
\medskip 

\noindent \textbf{Lemma 5.} \textit{Let $\rho _{B}:B\rightarrow B/H\subseteq 
\mathbb{R}^{\ell}$ be the orbit mapping of a linear action of a compact Lie group $H$ on an open ball $B$ in $\mathbb{R}^{n}.$ For every smooth vector field 
$Y$ on $\Sigma =B/H$, which extends to a derivation on $\overline{\Sigma}=\overline{B}/H$, there is a bounded smooth $H$ invariant vector field $X$ on $B$, which is $\rho _{B}$ related to $Y$ on $B$. In other words, 
$Y(\overline {f})=((\rho_{B})_{\ast }\comp X\comp (\rho _{B})^{\ast })(\overline{f})$ for every $\overline{f} \in C^{\infty }(B/H)$.}\medskip

\noindent \textbf{Proof.} Let $Y$ be a smooth vector field on $B/H$, which extends to a derivation on $\overline{\Sigma}=\overline{B}/H$. Since 
$\overline{B}/H$ is a closed differential subspace of $\mathbb{R}^{\ell}$, in coordinates $\mathbf{y}=(y_{1},...,y_{\ell})$ on $\mathbb{R}^{\ell}$, we may write 
\[
Y(\mathbf{y})=\sum_{i=1}^{\ell}g_{i}(\mathbf{y})\frac{\partial }{\partial y_{i}} , 
\]%
where $g_{i}\in C^{\infty }(B/H)$ is the restriction to $B/H$ of a smooth
function to $\mathbb{R}^{\ell}$. \medskip

We begin the proof by showing that the orbit space $B/H$ is connected.
Observe that the open ball $B\subseteq \mathbb{R}^{n}$ is centered at the
origin and the action of $H$ on $B$ is the restriction of the linear action
of $H$ on $T_{m}M,$ see Equation (\ref{16}). The linearity of the action of 
$H$ on $B$ implies that it commutes with scalar multiplication. Moreover, the
origin $\mathbf{0}\in \mathbb{R}^{n}$ is $H$ invariant so that it is an
orbit of $H$. Hence $\mathbf{\overline{0}}=\rho _{B}(\mathbf{0})\in B/H$. Let 
$\mathbf{x}\in B.$ For each $t\in \lbrack 0,1]$ and every $h\in H$, we have 
$h\cdot (t\mathbf{x})=t(h\cdot \mathbf{x}).$ Therefore, the line segment 
$[0,1]\rightarrow B:t\mapsto t\mathbf{x.}$ Thus the $H$ orbit through $%
\mathbf{0}$ and the $H$ orbit through $\mathbf{x}$ belong to the same
connected component of $B/H$. This implies that $B/H$ is connected.\medskip

The connectedness of $B/H$ ensures that there is a unique principal type 
$(K) $ whose corresponding orbit type $B_{(K)}$ is open and dense in $B$, see
Duistermaat and Kolk \cite[p.118]{duistermaat-kolk}. Moreover, the orbit space $B_{(K)}/H$ is a connected smooth manifold, and 
\[
\rho _{B_{(K)}}:B_{(K)}\rightarrow B_{(K)}/H:\mathbf{x}\mapsto H\cdot 
\mathbf{x} 
\]%
is a locally trivial fibration, whose fibre over 
$\rho _{\mid B_{(K)}}(\mathbf{x})$ is the $H$ orbit $H\cdot \mathbf{x}.$ Hence, 
for every $\mathbf{y}\in B_{(K)}/H\subseteq B/H$, there exists an open neighbourhood $V$ of $\mathbf{y}\in B_{(K)}/H$ such that $W=\rho ^{-1}(V)$ is trivial. In other words, there is a diffeomorphism $\tau :W\rightarrow H\times V$ such that $\rho _{\mid W}=\mathrm{pr}_{2}\comp \tau ,$ where 
$\mathrm{pr}_{2}:H\times V\rightarrow V$ is the projection on the second factor. This implies that there is a smooth $H$ invariant vector field $X_{W}$ on $W\subseteq B_{\left( K\right) },$ which is $\rho _{\mid W}$ related to the restriction $Y_{\mid V}$ of $Y$ to $V=\rho (W)$.\medskip

Repeating the above argument at each point $y\in B_{(K)}/H$ leads to a
covering $\{W_{\alpha }\}_{\alpha \in I}$ of $B_{(K)}$ by $H$ invariant open
subsets $W_{\alpha }$ of $B_{(K)}$ on which there exists an $H$ invariant
vector field $X_{W_{\alpha }},$ which is $\rho $ related to the restriction
of the vector field $Y$ to $V_{\alpha }=W_{\alpha }/H.$ Using an $H$
invariant partition of unity on $B_{(K)},$ we obtain a vector field 
$X_{B_{(K)}}$ on $B_{(K)},$ which is $\rho _{\mid B_{(K)}}$ related to 
$Y_{\mid B_{(K)}/H}$, i.e., $Y_{\mid B_{(K)}/H}=$ 
$(\rho _{\mid B_{(K)}})_{\ast }\comp X_{B_{(K)}}\comp 
(\rho _{\mid B_{(K)}})^{\ast }.$
\medskip

The module $\mathfrak{X}(B)^{H}$ of $H$ invariant smooth vector fields on $B$
is finitely generated by polynomial vector fields, see \cite{cushman}, and we denote a generating set by $\{X_{j}\}_{j=1}^{N}$. Hence, every $H$
invariant smooth vector field $X_{B}$ on $B$ is of the form 
$X_{B}=\sum_{j=1}^{N}h_{j}X_{j}$ for some $h_{1},\ldots ,h_{n}\in C^{\infty }(B)$. 
Similarly, every $K$ invariant smooth vector field on $B_{(K)}$ can be
written as $\sum_{j=1}^{N}f_{j}X_{j}$, where $f_{j}\in C^{\infty }(B_{(K)})$. Since $B_{(K)}$ is open and dense in $B$, a generic 
$f\in C^{\infty}(B_{(K)})$ need not extend to a smooth function on $B$. Therefore a generic vector field on $B_{(K)}$ need not extend to a smooth vector field on $B$. On the other hand, the $H$ invariant vector field $X_{B_{(K)}}$ on 
$B_{(K)},$ is obtained above from a smooth bounded vector field $Y$ on $B/H.$ Therefore 
\begin{equation}
X_{B_{(K)}}(\mathbf{x})=\sum_{j=1}^{N}k_{j\mid B_{(K)}}(\mathbf{x})X_{j}(%
\mathbf{x}),  
\label{19}
\end{equation}%
for each $\mathbf{x}\in B_{(K)}\subseteq B\subseteq \mathbb{R}^{n},$ where
every $k_{j}\in C^{\infty }(\mathbb{R}^{n})$, and $k_{j\mid B_{(K)}}$ is the
restriction of $k_{j}$ to $B_{(K)}.$\medskip

Since $B_{(K)}$ is open and dense in $B$, we may define 
\begin{equation}
X(\mathbf{x} )=\left\{ 
\begin{array}{cl}
X_{B_{(K)}} ( \mathbf{x} ) & \mbox{if $\mathbf{x} \in B_{(K)}$} \\
\sum^N_{i=1} \lim_{k \rightarrow \infty }k_{j} ( {\mathbf{x}}_k ) X_j 
( {\mathbf{x}}_{k} ), &  
\parbox[t]{1.75in}{where $\{ {\mathbf{x}}_k \} \subseteq B_{(K)}$ and \\
$\mathbf{x}= \lim_{k \rightarrow \infty} {\mathbf{x}}_k \in 
B \setminus B_{(K)}$,}
\end{array} \right.  
\label{20}
\end{equation}
provided that $\lim_{k \rightarrow \infty}k_{j}( {\mathbf{x}}_{k})$ exists
and is unique. Since the vector fields $X_{j}({\mathbf{x}}_k)$ are
smooth on $B$, 
\begin{displaymath}
\lim_{k\rightarrow \infty }k_{j}({\mathbf{x}}_{k})X_{j}({\mathbf{x}}_k)
=(\lim_{k\rightarrow \infty }k_{j}({\mathbf{x}}_k))X_{j}(\lim_{k\rightarrow \infty }
{\mathbf{x}}_k)\, \, \mbox{for every $j=1,\ldots ,N$.}
\end{displaymath}

Moreover, since $B_{(K)}$ is open and dense in $B,$ it is open and dense in 
$\overline{B}$, the closure of $B$. In Equation (\ref{19}), each function 
$k_{j\mid B_{(K)}}$ is the restriction to $B_{(K)}\subseteq B\subseteq 
\mathbb{R}^{n}$ of a smooth function $k_{j}$ on $\mathbb{R}^{n}.$ Moreover,
the choice of polynomial basis $\{X_{j}\}_{j=1}^{N}$ ensures that the right
hand side of Equation (\ref{19}) extends to the closure $\overline{B}$ of $B$.
Hence all the the limits $\lim_{k\rightarrow \infty }k_{j}
(\boldsymbol{x}_{k})$ in Equation (20) exist, and $X(\mathbf{x})$ is defined for all $\mathbf{x}\in B$.\medskip

We need to show that this definition of $X$ on $B$ depends only on 
$X_{B_{(K)}}$. Since each $k_{j}$ is continuous on $B$ and its first
partial derivatives are bounded on $\overline{B}$, it follows that $k_{j}$ are
uniformly continuous on $\overline{B}.\,$In particular, if $c:[0,1]\rightarrow 
\overline{B}$ is a smooth curve, such that $c([0,1))\subseteq B_{(K)}$ and 
$c(1)\in B$, then
\[
k_{j\mid B}(c(1))=k_{j\mid B}(c(0))+\int_{0}^{1}\frac{\partial k_{j\mid B}}{%
\partial t}(c(t))dt=k_{j\mid B_{(K)}}(c(0))+\int_{0}^{1}\frac{\partial
k_{j\mid B_{(K)}}}{\partial t}(c(t))dt. 
\]%
Thus, the values of $k_{j}$ on $B$ are uniquely determined by $k_{j\mid
B_{(K)}}$. Repeating this argument for all the first order partial
derivatives of $k_{j}$, we deduce that the first order partial derivatives
of $k_{j}$ on $B$ are uniquely determined by $k_{j\mid B_{(K)}}$ and its
first partial derivatives. Continuing this process for every partial
derivative of every order shows that the restriction of $k_{j}$ to $B$ is
uniquely determined by $k_{j\mid B_{(K)}}$.\medskip

The above argument applies to each of the functions $k_{i},$ for $i=1,\ldots ,n,$
in Equation (\ref{19}) and ensures that the $H$ invariant vector field 
$X_{B_{(K)}}$, thought of as the smooth section $B_{(K)}\rightarrow
TB_{(K)}=B_{(K)}\times \mathbb{R}^{n}$, extends to a smooth $H$ invariant
map $X:B\rightarrow B\times \mathbb{R}^{n}$, which is $\rho $ related to the
section $Y:B/H\rightarrow T(B/H).$\medskip

It remains to show that $X$ is a vector field on $B$. By construction, 
$X_{B_{(K)}}$ is an $H$ invariant vector field on an open dense subset 
$B_{(K)}$ of $B$, which is $\rho _{\mid B_{(K)}}$ related to the vector field 
$Y_{\mid B_{(K)}/H}$. The closure of $B_{(K)}$ in $B$ is the union of orbit
types $B_{(J)},$ where $(J)\leq (K).$ Suppose that $\boldsymbol{x} =\lim_{k\rightarrow \infty }\boldsymbol{x}_{k}\in B_{(J)},$ where 
$\boldsymbol{x}_{k}\in B_{(K)}$ for all $k\in \mathbb{Z}_{\geq 1}.$ \ Then 
$\boldsymbol{y}\in \rho _{B}(\boldsymbol{x})\in B_{(J)}/H$ and 
\begin{align*}
T_{\boldsymbol{x}}\rho _{B}(X(\boldsymbol{x})) &= T_{\boldsymbol{x}}
\rho _{B}(X(\lim_{k\rightarrow \infty }\boldsymbol{x}_{k}))=\lim_{k\rightarrow
\infty }T_{\boldsymbol{x}}\rho _{B}(X(\boldsymbol{x}_{k})) \\
&= \lim_{k\rightarrow \infty }Y(\rho _{B}(\boldsymbol{x}_{k}))=
Y(\rho _{B}(\boldsymbol{x}))=Y(\boldsymbol{y})\text{,}
\end{align*}%
because $Y_{\mid B_{(K)/H}}$ is the restriction to $B_{(K)}/H$ of a smooth,
and hence continuous, vector field on $B/H$. By Proposition 7, for every
orbit type $B_{(J)}$, the manifold $B_{(J)}/H$ is an invariant manifold of
of the vector field $Y$. So $Y_{\mid B_{(J)}/H}$ is a vector field on 
$B_{(J)/H}.$ Hence, for every $\boldsymbol{x}\in B$, 
\[
X(\boldsymbol{x})\in (T_{\boldsymbol{x}}\rho _{B})^{-1}
(Y(\boldsymbol{y}))\subseteq (T_{\boldsymbol{x}}\rho _{B})^{-1}(T_{\boldsymbol{y}}(B_{(J)}/H))\subseteq T_{\boldsymbol{x}}B.
\]
Therefore, $X$ is a smooth vector field on $B$, which is $\rho _{B}$ related
to the vector field $Y$ on $B/H$. \hfill $\square $ \medskip 

The aim of the rest of this section is to prove \medskip 

\noindent \textbf{Theorem 2.} \textit{Let $\Phi : G \times M \rightarrow M$ be a proper action of a Lie group $G$ on a connected smooth manifold $M$ with orbit map $\pi : M \rightarrow M/G$. Every smooth vector field on 
the locally closed subcartesian differential space $(M/G, C^{\infty}(M/G))$ is 
$\pi $ related to a smooth $G$ invariant vector field on $M$.} \medskip 

First we prove \medskip 

\noindent \textbf{Lemma 6.} \textit{Let $S_m$ be a slice to the $G$ action $\Phi $ at $m \in M$ and suppose that $\widehat{X}$ is a smooth $H$ invariant vector field on some $H$ invariant open neighborhood $U_m$ of $m$ in $S_m$. Here $H$ is the isotropy group $G_m$ at $m$. Then the vector field $\widehat{X}$ extends to a smooth $G$ invariant vector field $X$ on $M$.} \medskip 

\noindent \textbf{Proof.} Let $U_m \subseteq S_m$ be an $H$ invariant open subset of $S_m$ containing $m$. Because $S_m$ is a slice, $G\cdot U_m = \{ {\Phi }_g(U_m) \in M \setrule \, g \in G \}$ is a $G$ invariant open subset of $M$, which contains the $G$ orbit $G\cdot m$. On $G \cdot U_m$ define the vector field $X  = \{ ({\Phi }_g)_{\ast }\widehat{X} \setrule \, g \in G\}$. We check that $X$ is well defined. Suppose that $g \cdot s_m = g' \cdot s'_m$, where $g$, $g' \in G$ and $s_m$, $s'_m \in S_m$. Since $S_m$ is a slice, it follows that $g^{-1}g' = h \in H$. Hence 
\begin{displaymath}
({\Phi }_{g'})_{\ast }\widehat{X} = ({\Phi }_{gh})_{\ast }\widehat{X} = 
({\Phi }_g)_{\ast }\big( ({\Phi }_h)_{\ast } \widehat{X} \big) = ({\Phi }_g)_{\ast }\widehat{X}, 
\end{displaymath}
where the last equality above follows because the vector field $\widehat{X}$ is $H$ invariant. So the vector field $X$ on $G \cdot U_m$ is well defined and by definition is $G$ invariant.  
\par Next we show that $X$ is smooth. Let $L$ be a complement to the Lie algebra $\mathfrak{h}$ of the Lie group $H$ in the Lie algebra $\mathfrak{g}$ of the Lie group $G$. For every $\xi \in L$ and $\eta \in \mathfrak{h}$ consider the map $\mu : G \rightarrow (\exp L)H: \exp (\xi + \eta ) \mapsto \exp \xi \, 
\exp \eta $, which sends the identity element $e_G$ of $G$ to $e_G \cdot e_H = e_G$. It is a local diffeomorphism, since its tangent $T_{e_G}\mu : \mathfrak{g} \rightarrow L \oplus \mathfrak{h} = \mathfrak{g}$ is 
the identity map. Thus there are open subsets $V_G$, $V_L$, and $V_H$ of $e_G$, $0_L$, and $e_H$, respectively, such that $\mu (V_G) = \exp V_L \cdot V_H$. Hence every $g \in V_G$ may be written uniquely as $g = (\exp \xi )h$ for some $\xi \in V_L$ and some $h \in V_H$. For every $s_m \in U_m \subseteq S_m$ we have 
\begin{align}
({\Phi }_g)_{\ast }X(s_m) & = ({\Phi }_{(\exp \xi )h})_{\ast } X(s_m) = 
( {\Phi }_{\exp \xi })_{\ast } \big( ({\Phi }_h)_{\ast } X(s_m) \big) \notag \\
& = ({\Phi }_{\exp \xi })_{\ast} X(s_m), \, \, \mbox{since $X$ is $H$ invariant on $S_m$} \notag \\
& = T_{s_m}{\Phi }_{s_m} X \big( {\Phi }_{\exp -\xi }(s_m) \big) .  
\label{eq-s3one} 
\end{align}

Consider the local diffeomorphism 
\begin{displaymath}
\varphi : V_L \times S_m \rightarrow G \cdot S_m: (\xi , s) \mapsto 
{\Phi }_{\exp \xi }(s) = \Phi (\exp \xi , s). 
\end{displaymath}
Then ${\Phi }_{\exp \xi } = {\varphi }_{|\{ \xi \} \times S_m}$. So 
if $q = \varphi (\xi ,s) ={\Phi }_{\exp  \xi }(s)$, then 
${\Phi }^{-1}_{\exp \xi }(q) = s$. Let $W$ be a neighborhood of 
$\{ 0 \} \times S_m \subseteq V_L \times S_m$ such that $\varphi $ restricted to 
$W$ yields a diffeomorphism ${\varphi }_{|W}: W \subseteq V_L \times S_m 
\rightarrow U = \varphi (W) \subseteq G \cdot S_m$. For $s \in S_m$ let 
${\mathrm{e}}^{t {\widehat{X}}_m}(s)$ be the integral curve of the vector 
field ${\widehat{X}}_m$ starting at $s$. Since $X_m$ is a $G$ invariant extension 
of ${\widehat{X}}_m$ to a vector field on $G\cdot S_m$ (whose smoothness 
we want to prove) of a smooth $H$ invariant vector field ${\widehat{X}}_m$ on 
$S_m$, it follows that $X_m|S_m = {\widehat{X}}_m$ is a smooth vector field on $S_m$. Therefore 
\begin{equation}
{\varphi }^{{\widehat{X}}_m}_t(s) = {\varphi }^{X_m|S_m}_t(s) = 
{\varphi }^{X_m}_t(s), 
\label{eq-jedfive}
\end{equation}
for all $s \in S_m$. Consider a curve $c_q$ in $U \subseteq G \cdot S_m$ starting at $q = {\Phi }_{\exp \xi }(s)$ defined by 
$c_q(t) = {\Phi }_{\exp \xi }\big( {\varphi }^{{\widehat{X}}_m}_t 
({\Phi }^{-1}_{\exp \xi }(q)) \big) $. Using Equation (\ref{eq-jedfive}) we obtain 
${\dot{c}}_q(0) = T_q{\Phi }_{\exp \xi }({\widehat{X}}_m(s)) = X(q)$ for 
all $q \in U$. Since the family of curves $t \rightarrow c_q(t)$ depends smoothly 
on $q \in U$ and $U$ is an open subset of $G \cdot S_m$ containing $S_m$, it follows that $X_U$ is a smooth vector field on $U$. For any $m' \in G\cdot S_m$ there exists a $g \in G$ such that the open set ${\Phi }_g(U)$ contains $m'$. Since $X$ is $G$ invariant, smoothness of $X$ on $U$ ensures that $X$ is smooth on ${\Phi }_g(U)$. Hence $X$ is a smooth vector field on $G \cdot S_m$. \medskip 

The above argument can be repeated at each point $m \in M$. This leads to a covering ${\{ G \cdot S_{m_{\alpha}} \} }_{\alpha \in I}$ of $M$ by open $G$ invariant subsets $G \cdot S_{m_{\alpha }}$, where $S_{m_{\alpha }}$ is a slice at $m_{\alpha }$ for the action of $G$ on $M$ and $I$ is an index set. If $Y$ is a vector field on $\overline{M}$, then for each $\alpha \in I$ there exists a $G$ invariant vector field $X_{m_{\alpha}}$ on $G \cdot S_{m_{\alpha}}$ that is $\pi $ related to the restriction of $Y$ to $(G \cdot S_{m_{\alpha }})/G \subseteq M/G$. Using a $G$ invariant partition of unity on $M$ subordinate to the covering 
${\{ G \cdot S_{m_{\alpha }} \} }_{\alpha \in I}$, we can glue the pieces 
$X_{m_{\alpha }}$ together to obtain a smooth $G$ invariant vector field $X$ on $M$, which is $\pi $ related to the vector field $Y$ on $M/G$. \hfill $\square $ \medskip 

\noindent \textbf{Proof of theorem 2.} Applying Lemma 6 to the push forward 
by the local diffeomorphism ${\psi }_{|B}: B \subseteq T_mM \rightarrow U_m \subseteq S_m \subseteq M$, given by the Bochner lemma, of the vector field on $B$ constructed in Lemma 5, proves Theorem 2.  \hfill $\square $ \medskip 

\noindent \textbf{Proposition 9} \textit{If $Y$ is a derivation of $C^{\infty}(M/G)$, which is $\pi $ related to a derivation $X$ of $C^{\infty}(M)^G$, then $Y$ is a smooth vector field on $M/G$.} \medskip 

\noindent \textbf{Proof.} Since $M$ is a smooth manifold, $X$ is a smooth $G$ invariant vector field on $M$, which is $\pi $ related to derivation $Y$ of $C^{\infty}(M/G)$. Thus the image under 
$\pi $ of a maximal integral curve of $X$ on $M$, is a maximal integral curve of $Y$ on $M/G$. Hence $Y$ is a smooth vector field on the locally closed subcartesian differential space $\big( M/G , C^{\infty}(M/G) \big) $. 
\hfill $\square $  

\section{Differential $\mathbf{1}$-forms on the orbit space}

In this section we define the notion of a differential $1$-form on 
the orbit space $M/G$ of a proper group action $\Phi : G \times M \rightarrow M: 
(m,g) \mapsto g \cdot m$ on a smooth manifold $M$ with orbit map 
$\pi : M \rightarrow M/G: m \mapsto \overline{m} = G\cdot m$.  We show that the differential $1$-forms on $M/G$ together with the exterior derivative generate a differential exterior algebra.  \medskip 

Theorem 2 and Proposition 9 show that $Y$ is a vector field on 
$M/G$ if and only if there is a $G$ invariant vector field $X$ on $M$, which is 
$\pi $ related to $Y$, that is, every integral curve of $Y$ is the image under 
the map $\pi $ of an integral curve of $X$. Let ${\Lambda }^1(M/G)$ 
be the set of differential $1$-forms on $M/G$, that is, the set of linear mappings
\begin{displaymath}
\theta : \mathfrak{X}(M/G) \rightarrow C^{\infty}(M/G): Y \mapsto 
\theta (Y) = \langle \theta | Y \rangle , 
\end{displaymath} 
which are linear over the ring $C^{\infty}(M/G)$, that is, 
$\langle \theta | \overline{f}\, Y \rangle = \overline{f} \, \langle  \theta | Y \rangle $ 
for every $\overline{f} \in C^{\infty}(M/G)$ and every 
$Y \in \mathfrak{X}(M/G)$. \medskip 

In order to prove some basic properties of differential $1$-forms on $M/G$, we 
need to prove some properties of the $G$ orbit map $\pi $ (\ref{eq-s0zeronw}). 
\medskip 

The map 
\begin{displaymath}
\begin{array}{l}
T_m\pi : T_m M \rightarrow T_{\overline{m}}(M/G)= {\spann }_{\R }\{ Y(\overline{m}) 
\setrule \, Y \in \mathfrak{X}(M/G) \} : \\
\rule{0pt}{12pt}\hspace{.5in} v_m = X(m) \mapsto Y(\overline{m}),
\end{array} 
\end{displaymath}
where $X \in \mathfrak{X}(M)^G$ and $Y$ is the vector field on $M/G$ constructed in Proposition 2, is the \emph{tangent} of the map $\pi $ at $m \in M$. To show that $T_m \pi $ is well defined we argue as follows. Suppose that $v_m = X'(m)$, where $X' \in \mathfrak{X}(M)^G$. Then 
\begin{displaymath}
T_m\pi \big( X(m) - X'(m) \big) = T_m \pi (v_m -v_m) = 0, 
\end{displaymath}
since $T_m \pi $ is a linear map. \hfill $\square $ \medskip 

\noindent \textbf{Lemma 7.} \textit{For each $m \in M$ 
\begin{equation}
\ker T_m \pi = {\spann }_{\R } \{ X_{\xi }(m) \in T_mM \setrule \, \xi \in \mathfrak{g} \} ,
\label{eq-s4fourvnw}
\end{equation}
where $\mathfrak{g}$ is the Lie algebra of $G$.} \medskip 

\noindent \textbf{Proof.} By definition ${\pi }^{-1}(\overline{m}) = G \cdot m$. Thus 
\begin{equation}
T_m\big( {\pi }^{-1}(\overline{m}) \big) = 
T_m(G \cdot m) = {\spann}_{\R } \{ X_{\xi} \in T_mM \setrule \, \xi \in \mathfrak{g} \} . 
\label{eq-s4fourvnwstar}
\end{equation}
The curve ${\gamma }_m: I_m \subseteq \R \rightarrow M: t \mapsto \exp t\xi \cdot m$ is an integral curve of $X_{\xi }$ starting at $m$. So 
$\pi \big( {\gamma }_m(t) \big) = 
\pi (m) = \overline{m}$ for every $t \in I_m$. Thus $T_m \pi X_{\xi }(m) = 
\smalldbydt  \hspace{-10pt} \pi \big( {\gamma }_m(t) \big) = 0_{\overline{m}}$, that is, $X_{\xi }(m) \in \ker T_m\pi $. Consequently, 
${\spann}_{\R } \{ X_{\xi} \in T_mM \setrule \, \xi \in \mathfrak{g} \} \subseteq \ker T_mM$. 
\par To prove the reverse inclusion, we argue as follows. Since $m \in M$, it follows that 
$m \in M_{(H)}$, where $H = G_m$. Let $M_{(K)}$ be the maximal orbit type of the proper $G$ action on $M$. The maximal orbit type $M_{(K)}$ is a dense open subset of $M$, whose boundary 
$\partial M_{(K)} = \mathrm{cl}(M_{(K)}) \setminus M_{(K)}$ contains $M_{(H)}$,  
since the orbit types of the $G$ action stratify $M$. Suppose that $v_m$ is a nonzero vector in $\ker T_m\pi $. There is a vector field $X$ on $M$ with $X(m) = v_m$ having an integral curve ${\gamma }_p:I_p \subseteq \R \rightarrow M: t \mapsto {\varphi }^X_t(p)$ 
starting at $p \in M_{(K)}$ such that ${\gamma }_p(\tau ) = m$ for some $\tau \in 
I_p \cap {\R }_{>0}$. We may suppose that ${\gamma }_p([0, \tau )) \subseteq M_{(K)}$. Since ${\overline{M}}_{(K)} = \pi (M_{(K)})$ is a smooth submanifold of the differential 
space $(M/G, C^{\infty}(M/G))$, the curve 
${\Gamma }_m:(0 , \tau ] \rightarrow {\overline{M}}_{(K)}: t \mapsto  
{\gamma }_p(\tau -t)$ is a smooth integral curve of the vector field $-X$ such that 
${\Gamma }_m((0, \tau ]) \subseteq M_{(K)}$. Hence on $(0, \tau ]$ the curve $\pi \comp {\Gamma }$ on the smooth manifold ${\overline{M}}_{(K)}$ is smooth. Thus 
\begin{align*}
\frac{\dee }{\dee t}(\pi \comp {\Gamma }_m)(t) & = 
\frac{\dee }{\dee t} \pi \big( {\varphi }^{-X}_t(m) \big) = 
\dbyds \pi \big( {\varphi }^{-X}_{t+s}(m) \big) 
= \dbyds \pi \big( {\varphi }^{-X}_t({\varphi }^{-X}_s(m)) \big) \\
& = T_{{\varphi }^{-X}_t(m)}\pi \big( T_m\pi (-X(m)) \big) 
= 0_{\pi ({\varphi }^{-X}_t(m))},  
\end{align*} 
since $X(m) = v_m \in \ker T_m \pi $. Thus the curve $\pi \comp {\Gamma }_m$ is constant, since ${\overline{M}}_{(K)}$ is a smooth manifold. Because the curve $\pi \comp {\Gamma }_m$ 
is continuous on $[0, \tau ]$, we get $\pi \big( {\Gamma }_m(t)\big) = \pi \big( {\Gamma }_m(0) \big) 
= \pi (m) = \overline{m}$. Hence ${\Gamma }_m(t) \subseteq {\pi }^{-1}(\overline{m})$ for 
all $t \in [0, \tau ]$. But $\lim_{t \searrow 0}{\Gamma }^{\prime}_m(t) = v_m$. So 
$v_m \in T_m{\pi }^{-1}(\overline{m})$. Hence 
\begin{displaymath}
\ker T_m\pi \subseteq T_m{\pi }^{-1}(\overline{m}) = 
{\spann}_{\R }\{ X_{\xi }(m) \in T_mM \setrule \, \xi \in \mathfrak{g} \} ,
\end{displaymath}
where the equality follows from Equation (\ref{eq-s4fourvnwstar}). This verifies Equation (\ref{eq-s4fourvnw}). \hfill $\square $ \medskip 

A differential $1$-form $\omega $ on $M$ is \emph{semi-basic} with respect to 
the $G$ action $\Phi $ if and only if $X_{\xi } \lefthook \, \omega =0$ for 
every $\xi \in \mathfrak{g}$, the Lie algebra of $G$. \medskip 

\noindent \textbf{Proposition 10.} \textit{For every $\theta \in {\Lambda }^1(M/G)$, the differential $1$-form ${\pi }^{\ast } \theta $ on $M$ is $G$ invariant and semi-basic.} \medskip 

\noindent \textbf{Proof.} By definition of ${\pi }^{\ast }$ the map 
\begin{equation}
{\pi }^{\ast } \theta : \mathfrak{X}(M)^G \rightarrow C^{\infty}(M)^G: 
X \mapsto \langle {\pi }^{\ast }\theta | X \rangle 
\label{eq-s4vnwone}
\end{equation}
is linear, since the map 
$\mathfrak{X}(M)^G \rightarrow \mathfrak{X}(M/G): X \mapsto Y$ 
is linear. Moreover, for any $f \in C^{\infty}(M)^G$ 
\begin{align*}
\langle {\pi }^{\ast }\theta | fX \rangle & = 
{\pi }^{\ast }(\langle \theta | \overline{f}\, Y \rangle ), \, \, 
\parbox[t]{3in}{since the map $X \mapsto Y$ is a module \\ homomorphism} \\
& = {\pi }^{\ast }(\overline{f}\, \langle \theta | Y \rangle ), \, \, 
\mbox{because $\theta \in {\Lambda }^1(M/G)$} \\
& = {\pi }^{\ast }(\overline{f}) \, {\pi }^{\ast }(\langle \theta | Y \rangle ) = 
f \, \langle {\pi }^{\ast }\theta | X \rangle . 
\end{align*}
Thus ${\pi }^{\ast }\theta \in {\Lambda }^1(M/G)$. For every $\xi \in 
\mathfrak{g}$ 
\begin{displaymath}
L_{X_{\xi }}( \langle {\pi }^{\ast }\theta | X \rangle ) = 
L_{X_{\xi }}\big( {\pi }^{\ast }( \langle \theta | Y \rangle ) \big) = 0, 
\end{displaymath}
because ${\pi }^{\ast }( \langle \theta | Y \rangle \in C^{\infty}(M)^G$. So 
${\pi }^{\ast }\theta $ is a semi-basic $1$-form on $M$. \hfill $\square $ \medskip 

\noindent \textbf{Proposition 11.} \textit{Let $\vartheta $ be a $G$ invariant semi-basic differential $1$-form on $M$. Then there is a $1$-form $\theta $ on $M/G$ 
such that $\vartheta = {\pi }^{\ast }\theta $.} \medskip 

\noindent \textbf{Proof.} Given $Y \in \mathfrak{X}(M/G)$, there is 
an $X \in \mathfrak{X}(M)^G$, which is $\pi $ related to $Y$, that is, 
$T_m \pi X(m) = Y(\pi (m))$ for every $m \in M$. It is clear that the 
definition of $\theta $ needs to be 
\begin{equation}
{\pi }^{\ast }(\langle \theta | Y \rangle ) = \langle \vartheta | X \rangle . 
\label{eq-s4twonw}
\end{equation}
It remains to show that $\theta $ is well defined. Since the $1$-form $\vartheta $ and the vector field $X$ are $G$ invariant, we get 
\begin{displaymath}
{\Phi }^{\ast}_g\big( \langle \vartheta | X \rangle \big) (m)  = 
\langle {\Phi }^{\ast }_g\vartheta | {\Phi }^{\ast }_gX \rangle ({\Phi }_g(m)) = 
\langle \vartheta | X \rangle (m), 
\end{displaymath}
for every $(g,m) \in G \times M$. Thus the function $M \rightarrow \R : m \mapsto 
\langle \vartheta | X \rangle (m)$ is smooth and $G$ invariant. We now show that the mapping 
$\theta : \mathfrak{X}(M/G) \rightarrow C^{\infty}(M/G)$, 
where $\theta $ is given in Equation (\ref{eq-s4twonw}), is well defined. Suppose 
that $X' \in \mathfrak{X}(M)^G$ such that $X'$ is $\pi $ related to $Y$. Then 
$T_m\pi (X(m) - X'(m)) = Y(m) - Y(m) =0$ for every $m \in M$. So 
$(X(m) - X'(m)) \in {\spann }_{\R }\{ X_{\xi }(m) \in T_mM \setrule \, 
\xi \in \mathfrak{g} \} $, by Proposition 10. Thus 
\begin{align*}
\langle \vartheta | X \rangle & = \langle \vartheta | (X - X') \rangle + 
\langle \vartheta | X' \rangle = \langle \vartheta | X' \rangle , 
\end{align*}
since the $1$-form $\vartheta $ on $M$ is semi-basic. This shows that the 
map $\theta : \mathfrak{X}(M/G) \rightarrow C^{\infty}(M/G)$ is well defined. From Equation (\ref{eq-s4twonw}) it follows that $\theta $ is a linear mapping and that 
$\langle \theta | \overline{f}\, Y \rangle = \overline{f} \langle \theta | Y \rangle $ for 
every $\overline{f} \in C^{\infty}(M/G)$. Hence $\theta $ is a differential $1$-form 
on $M/G$, that is, $\theta \in {\Lambda }^1(M/G)$. Every 
$X \in \mathfrak{X}(M)^G$ is $\pi $ related to a $Y \in 
\mathfrak{X}(M/G)$. Thus $\langle {\pi }^{\ast }\theta | X \rangle  =  
{\pi }^{\ast }(\langle \theta | Y \rangle ) = \langle \vartheta | X \rangle $,  
that is, $\vartheta = {\pi }^{\ast }\theta $. \hfill $\square $  

\section{de Rham's theorem}

In this section we construct an exterior algebra of differential forms on 
the orbit space $M/G$ with an exterior derivative $\dee $ and 
show that de Rham's theorem holds for the sheaf of differential exterior algebras. 
\medskip 

Let $\ell \in {\Z }_{\ge 1}$. A differential $\ell $-form $\theta $ on $M/G$ is an element of  $L^{\ell }_{\mathrm{alt}}(T(M/G), \R )$, the set of 
alternating $\ell $ multilinear real valued mappings on $T(M/G)= 
\mathfrak{X}(M/G)$, namely,  
\vspace{-.1in}
\begin{displaymath}
\begin{array}{l}
\theta : \overbrace{\mathfrak{X}(M/G) \times \cdots \times \mathfrak{X}(M/G)}^{\ell } 
\rightarrow C^{\infty}(M/G): \\
\hspace{.65in} (Y_1, \ldots , Y_{\ell }) \longmapsto 
Y_{\ell }\lefthook ( \cdots \lefthook (Y_1 \lefthook \theta ) \ldots ) = 
\langle \theta | (Y_1, \ldots , Y_{\ell }) \rangle,
\end{array} 
\end{displaymath}
which is linear over $C^{\infty}(M/G)$, that is,  
\begin{displaymath}
\langle \theta | \big( Y_1, \ldots, \overline{f} \, 
Y_i, \ldots , Y_{\ell }) \rangle  = 
\overline{f} \langle  \theta | \big( Y_1, \ldots , 
Y_{\ell } \big) \rangle  
\end{displaymath} 
for every $1 \le i \le \ell $ and every 
$\overline{f} \in C^{\infty}(M/G)$. A differential $0$-form on $M/G$ is a smooth function on $M/G$. Let ${\Lambda }^{\ell }(M/G)$ be the real vector space of differential $\ell $-forms on $M/G$. For each $\overline{m} \in M/G$ let ${\Lambda }^{\ell }_{\overline{m}}(M/G) = 
{\spann }_{\R }\{ \theta (\overline{m}) \in L^{\ell }_{\mathrm{alt}}
(T_{\overline{m}}(M/G), \R ) \setrule \, \theta \in {\Lambda }^{\ell }(M/G) \} $. \medskip

\noindent \textbf{Proposition 12.} \textit{Let $\theta \in {\Lambda }^{\ell }(M/G)$ with $\ell \in {\Z }_{\ge 1}$. Then the  $\ell $-form $\vartheta = {\pi }^{\ast }\theta \in 
{\Lambda }^{\ell }_{\mathrm{sb}}(M)^G$, the set of semi-basic $G$ invariant 
$\ell $-forms on $M$. Here 
\begin{displaymath}
({\pi }^{\ast} \theta )(m) \big( X_1(m), \ldots , X_{\ell }(m) \big) = 
\theta (\pi (m)) \big( T_m \pi X_1(m), \ldots , T_m\pi X_{\ell }(m) \big) , 
\end{displaymath}
for every $m \in M$ and every $X_j \in \mathfrak{X}(M)^G$ for $1 \le j \le \ell $.} \medskip 

\noindent \textbf{Proof.} The proof is analogous to the proof of Proposition 10 for 
$1$-forms on $M/G$ and is omitted. \hfill $\square $ \medskip 

\noindent \textbf{Proposition 13.} \textit{Let 
$\vartheta \in {\Lambda }^{\ell }_{\mathrm{sb}}(M)^G$, 
where $\ell \in {\Z }_{\ge 1}$. Then there is an $\ell $-form $\theta \in 
{\Lambda }^{\ell}(M/G)$ such that $\vartheta = {\pi }^{\ast }(\theta )$.} \medskip 

\noindent \textbf{Proof.} The proof is analogous to the proof of Proposition 11 for 
$G$ invariant semi-basic $1$-forms on $M$ and is omitted. \hfill $\square $ \medskip 

We now define the exterior algebra $\Lambda (M/G)$ of differential forms on $M/G$. Let $\theta \in {\Lambda }^{h}(M/G)$ and $\phi \in {\Lambda}^k(M/G)$. The \emph{exterior product} is the $h+k$ form $\theta \wedge \phi $ on $M/G$ corresponding to the $G$ invariant semi-basic $h+k$-form ${\pi }^{\ast }\theta \wedge {\pi }^{\ast}\phi $ on $M$. Then 
$\big( \Lambda (M/G) = \sum_{\ell } \oplus {\Lambda }^{\ell }(M/G), \wedge \big) $ is an exterior algebra of differential forms on $M/G$. \medskip 

The exterior derivative operator $\dee $ on $\Lambda (M/G)$ is defined in terms 
of the Lie bracket of vector fields on $M/G$. If $Y$, $Y' \in 
\mathfrak{X}(M/G)$, then there are $X_Y$, ${X'}_{Y'} \in \mathfrak{X}(M)^G$, each of which is $\pi $ related to $Y$ and $Y'$, respectively. Their Lie bracket 
$[X_Y, {X'}_{Y'}] \in \mathfrak{X}(M)^G$. Then there is a vector field 
$Y_{[X_Y, {X'}_{Y'} ]}$ on $M/G$, which is $\pi $ related to $[X_Y, {X'}_{Y'}]$. 
Define $[Y, Y'] = Y_{[X_Y, {X'}_{Y'} ]}$. The following lemma shows that 
this Lie bracket is well defined. \medskip 

\noindent \textbf{Lemma 8.} \textit{For every $\overline{f} \in C^{\infty}(M/G)$ and 
every $Y$, $Y' \in \mathfrak{X}(M/G)$ 
\begin{equation}
[Y,Y'] (\overline{f}) = Y'(Y(\overline{f})) - Y(Y'(\overline{f})). 
\label{eq-s4nweightstar}
\end{equation} }

\noindent \textbf{Proof.} We compute. 
\begin{align*}
{\pi }^{\ast }\big( [Y,Y'](\overline{f}) \big) & = [X_Y, X_{Y'}]\big( {\pi }^{\ast }(\overline{f}) \big), \, \, 
\mbox{by definition of Lie bracket} \\
& \hspace{-.55in} = X_{Y'}\big( X_Y ( {\pi}^{\ast }(\overline{f})) \big) - 
X_Y\big( X_{Y'} ( {\pi}^{\ast }(\overline{f})) \big) , \, \, \mbox{because $X_Y$, $X_{Y'} \in 
\mathfrak{X}(M)^G$} \\
& \hspace{-.55in} = X_{Y'}\big( {\pi }^{\ast }(Y(\overline{f})) \big) - 
X_Y\big( {\pi }^{\ast }(Y'(\overline{f})) \big), \\
&\hspace{-.25in} \parbox[t]{3.6in}{since $X_Y({\pi }^{\ast }(\overline{f})) = {\pi }^{\ast }\big( Y(\overline{f}) \big)$ 
and $X_{Y'}({\pi }^{\ast }(\overline{f})) = {\pi }^{\ast }\big( Y'(\overline{f}) \big)$} \\
&\hspace{-.55in} = {\pi }^{\ast }\big( Y'(Y( \overline{f}) ) \big) - {\pi }^{\ast }\big( Y(Y'( \overline{f}) ) \big) \\
& \hspace{-.55in} = {\pi }^{\ast } \big( Y'(Y(\overline{f})) - Y(Y(\overline{f})) \big) ,  
\end{align*}
from which Equation (\ref{eq-s4nweightstar}) follows, because the orbit map 
$\pi $ is surjective. \hfill $\square $ \medskip 

\noindent \textbf{Corollary 2.} \textit{$[ \, \, . \, \, ]$ is a Lie bracket on 
$\mathfrak{X}(M/G)$.} 
\medskip 

\noindent \textbf{Proof.} The corollary follows from a computation using Equation 
(\ref{eq-s4nweightstar}). We give another argument. 
Bilinearity of the Lie bracket is straightforward to verify. We need only show that the Jacobi identity holds. We compute. 
\begin{align*}
[Y'', [Y, Y'] ] & = [Y_{X''}, Y_{[X,X']} ] = Y_{[ X_{Y''}, [ X_Y, X_{Y'}]]}\\
& = Y_{[ [X_{Y''}, X_Y], X_{Y'} ] +[X_Y, [X_{Y''}, X_{Y'}]]} \, \, 
\mbox{by the Jacobi identity on $\mathfrak{X}(M)^G$} \\
& = Y_{[ [X_{Y''}, X_Y], X_{Y'} ]} + Y_{[X_Y, [X_{Y''}, X_{Y'}]]} \\
& = [[Y'', Y], Y'] + [Y, [Y'', Y']], 
\end{align*}
which is the Jacobi identity on $\mathfrak{X}(M/G)$. \hfill $\square $ \medskip 

Let $\varphi $ be an $\ell $-form on $M/G$. Inductively define the exterior 
derivative $\dee $ of $\varphi $ as the $(\ell +1)$-form given by 
\begin{align}
\dee \varphi \big( Y_0, Y_1, \ldots , Y_{\ell } \big) & = 
\sum^{\ell }_{i=0} (-1)^i \dee \, (Y_i \lefthook \varphi ) \big( Y_0, \ldots , 
{\widehat{Y}}_i, \ldots , Y_{\ell } \big) \notag \\
& \hspace{-.25in} + \sum_{0\le i < j \le \ell }(-1)^{i+j}([Y_i, Y_j] \lefthook \varphi ) 
(Y_0, \ldots , {\widehat{Y}}_i, \ldots , {\widehat{Y}}_j, \ldots , Y_{\ell } ) .
\label{eq-s5one}
\end{align}
Here $Y_i \lefthook {\varphi }$ is the $\ell -1$ form on $M/G$ defined by 
\begin{displaymath}
\langle (Y_i \lefthook {\varphi }) |(Y_0, \ldots , {\widehat{Y}}_i, \ldots ,Y_{\ell -1}) \rangle = 
\langle \varphi | (Y_i, Y_1, \ldots , {\widehat{Y}}_i, \ldots , Y_{\ell -1}) \rangle , 
\end{displaymath} 
for $Y_0, \ldots , {\widehat{Y}}_i, \ldots , Y_{\ell -1} \in \mathfrak{X}(M/G)$. To complete the definition of exterior derivative, we define $\dee $ on $0$-forms. This we do as follows. Let $\overline{f} \in {\Lambda }^0(M/G) = C^{\infty}(M/G)$. Define 
the $1$-form $\dee \overline{f}$ by
\begin{displaymath}
\dee \overline{f}(\overline{m})\big( Y(\overline{m}) \big) = 
Y(\overline{f})(\overline{m}), 
\end{displaymath}
for every $\overline{m} \in M/G$, every $\overline{f}\in C^{\infty}(M/G)$, and every 
$Y \in \mathfrak{X}(M/G)$. \medskip 

\noindent \textbf{Lemma 9.}\textit{Let $\theta \in {\Lambda }^{\ell }(M/G)$. Then 
\begin{equation}
\dee \, ({\pi }^{\ast }\theta ) = {\pi }^{\ast }(\dee \theta ) .
\label{eq-s5onevnw}
\end{equation} }

\noindent \textbf{Proof.} Suppose that $\theta $ is an $\ell $-form on $M/G$. 
Pulling back the forms on both sides of Equation (\ref{eq-s5one}) by 
the orbit map $\pi $ gives  
\begin{align}
{\pi }^{\ast }(\dee \theta )(X_0, \ldots , X_{\ell }) & = 
\sum^{\ell }_{i=0}{\pi }^{\ast }(\dee \, (Y_i \lefthook \theta ))(X_0, \ldots , {\widehat{X}}_i, 
\ldots , X_{\ell }) \notag \\
& \hspace{-.75in} + \sum_{0\le i < j \le \ell }(-1)^{i+j}{\pi }^{\ast }([Y_i,Y_j] \lefthook \theta ) 
(X_0, \ldots , {\widehat{X}}_i, \ldots , {\widehat{X}}_j, \ldots , X_{\ell }) . 
\label{eq-s5onedagger}
\end{align}
By induction, assume that Equation (\ref{eq-s5onevnw}) holds for all 
forms of degree strictly less than $\ell $. Then 
\begin{equation}
{\pi }^{\ast }(\dee \, (Y_i \lefthook \theta )) = \dee \, ({\pi }^{\ast }(Y_i \lefthook \theta )). 
\label{eq-s5onestardagger}
\end{equation}
Now ${\pi }^{\ast }(Y_i \lefthook \theta ) = X_i \lefthook {\pi }^{\ast }\theta $, where 
$T\pi \comp X_i = Y_i \comp \pi$, since 
\begin{align*}
(X_i \lefthook {\pi }^{\ast }\theta )(X_0, \ldots , {\widehat{X}}_i, \ldots , X_{\ell }) & = 
(-1)^i ({\pi }^{\ast }\theta )(X_0, \ldots , X_i, \ldots , X_{\ell }) 
\\
&\hspace{-1.5in} = (-1)^i \theta (T\pi X_0, \ldots , T\pi X_i, \ldots , T\pi X_{\ell }) =
(-1)^i \theta (Y_0, \ldots ,Y_i, \ldots ,Y_{\ell }) 
\\
& \hspace{-1.5in} = (Y_i \lefthook \theta )(Y_0, \ldots , {\widehat{Y}}_i, \ldots , Y_{\ell }) 
={\pi }^{\ast }(Y_i \lefthook \theta )(X_0, \ldots , {\widehat{X}}_i, \ldots , X_{\ell }).  
\tag*{$\square $}
\end{align*}
Also 
\begin{align*}
{\pi }^{\ast }([Y_i,Y_j] \lefthook \theta )(X_0, \ldots , {\widehat{X}}_i, \ldots , 
{\widehat{X}}_j, \ldots , X_{\ell} ) & = \\
& \hspace{-1.5in} = ( [Y_i,Y_j] \lefthook \theta )(T\pi X_0, \ldots , 
{\widehat{T\pi X}}_i, \ldots , {\widehat{T\pi X}}_j , \ldots , T\pi X_{\ell } ) 
\\
& \hspace{-1.5in} = ([X_i, X_j] \lefthook {\pi }^{\ast } \theta )(X_0, \ldots , 
{\widehat{X}}_i, \ldots , {\widehat{X}}_j, \ldots , X_{\ell }), \\
& \hspace{-.75in}\mbox{since $[Y_i,Y_j] = [T\pi X_i, T\pi X_j] = 
T\pi [X_i, X_j]$.} 
\end{align*}
Thus Equation (\ref{eq-s5onedagger}) reads
\begin{align*}
{\pi }^{\ast }(\dee \theta )(X_0, \ldots , X_{\ell }) & = 
\sum^{\ell }_{i=0} (-1)^i \dee \, (X_i \lefthook {\pi }^{\ast }\theta ) 
(X_0, \ldots , {\widehat{X}}_i, \ldots , X_{\ell }) 
\\
& \hspace{-.5in}+ \sum_{0\le i < j \le \ell}(-1)^{i+j} ([X_i,X_j] \lefthook {\pi }^{\ast }\theta ) 
(X_0, \ldots, {\widehat{X}}_i, \ldots , {\widehat{X}}_j, \ldots ,X_{\ell }) 
\\
& \hspace{-.5in} = \dee \, ( {\pi }^{\ast }\theta )(X_0, \ldots , X_{\ell }). \tag*{$\square $}
\end{align*}

\noindent \textbf{Lemma 10.} \textit{If $\theta \in {\Lambda }^k(M/G)$ and 
$\phi \in {\Lambda }^h(M/G)$, then 
\begin{equation}
\dee \, (\theta \wedge \phi ) = \dee \theta \wedge \phi + (-1)^k \theta \wedge 
\dee \phi . 
\label{eq-s4onestarvnw}
\end{equation} }

\noindent \textbf{Proof.} On $M$ we have 
\begin{align*}
{\pi }^{\ast }(\dee \, (\theta \wedge \phi )) & = \dee \, ({\pi }^{\ast }(\theta \wedge \phi )) = 
\dee \, ({\pi }^{\ast }\theta \wedge {\pi }^{\ast }\phi ) \\
& = \dee \, ({\pi }^{\ast } \theta ) \wedge {\pi }^{\ast }\phi + 
(-1)^k {\pi }^{\ast }\theta \wedge \dee \, ({\pi }^{\ast }\phi ) \\
& = {\pi }^{\ast } \big( \dee \theta \wedge \phi + (-1)^k \theta \wedge \dee \phi \big) , 
\end{align*}
which implies that Equation (\ref{eq-s4onestarvnw}) holds, since 
the orbit map $\pi $ is surjective. \hfill $\square $ \medskip

\noindent \textbf{Lemma 11.} \textit{${\dee }^{\, 2}\theta  =0$ for every $\theta  \in 
{\Lambda }^{\ell }(M/G)$.} \medskip 

\noindent \textbf{Proof.} Suppose that $\ell \ge 1$. Then 
${\pi }^{\ast }\theta  $ is an $\ell $-form on $M$. Because $M$ is a 
smooth manifold, one has ${\dee }^{\, 2} ({\pi }^{\ast }\theta  ) =0$. 
By Lemma 9 ${\pi }^{\ast }(\dee \theta  ) = \dee \, ( {\pi }^{\ast }\theta  )$. So 
\begin{align*}
& {\pi }^{\ast }( {\dee }^{\, 2} \theta  ) = {\pi }^{\ast }(\dee \, (\dee \, \theta  )) 
= \dee \, ( {\pi }^{\ast }(\dee \theta  )) = {\dee}^{\, 2}({\pi }^{\ast }\theta  ) =0.
\end{align*}
Since the $G$ orbit map $\pi $ is surjective, ${\pi }^{\ast }({\dee }^{\, 2}\theta  ) =0$ implies ${\dee }^{\, 2}\theta  =0$. 
\par We now treat the case when $\ell =0$. Let $\overline{f} \in C^{\infty}(M/G)$ and let $Y_0$, $Y_1 \in \mathfrak{X}(M/G)$. Then 
\begin{align}
\dee \, ( \dee \overline{f})(Y_0,Y_1) & = 
\dee \, (Y_0 \lefthook \dee \overline{f})Y_1- \dee \, (Y_1 \lefthook \dee \overline{f})Y_0 
- [Y_0, Y_1]  \lefthook \dee \overline{f} \notag \\
& = \dee \, (Y_0(\overline{f}))Y_1 - \dee \, (Y_1(\overline{f}))Y_0 - 
[Y_0,Y_1](\overline{f}) \notag \\
& = Y_1(Y_0(\overline{f})) - Y_0(Y_1(\overline{f})) - Y_1(Y_0(\overline{f})) 
+ Y_0(Y_1(\overline{f})) =0. \tag*{$\square $}
\end{align}

We prove an equivariant version of the Poincar\'{e} lemma in ${\R }^n$.  \medskip

\noindent \textbf{Lemma 12.} \textit{Let $G$ be a Lie group, 
which acts linearly on ${\R }^n$ by $\Phi :G \times {\R }^n \rightarrow {\R }^n$. Let $H$ be a compact subgroup of $G$. Let $\beta $ be an $H$ invariant closed $\ell $-form with $\ell \ge 1$ on an open $H$ invariant ball $B$ centered at the origin of ${\R }^n$, whose closure is compact. Suppose that $\beta $ is semi-basic with respect to the $G$ action $\Phi $, that is, $X_{\xi } \lefthook \beta =0$ for every $\xi \in \mathfrak{g}$, the Lie algebra of $G$. Here $X_{\xi }(x) = T_e{\Phi }_m \xi$. Then there is an $H$ invariant $(\ell -1)$-form $\alpha $ on $B$, which is semi-basic with respect to the $G$ action $\Phi $, such that 
$\beta = \dee \alpha $.} \medskip 

\noindent \textbf{Proof.} Let $X$ be a linear vector field on ${\R }^n$ all of whose 
eigenvalues are negative real numbers. By averaging over the compact 
group $H$, we may assume that $X$ is $H$ invariant. Let ${\varphi }_t$ be the 
flow of $X$, which maps $B$ into itself. Moreover, ${\varphi }_{\infty} =0$. On $B$ one has 
\begin{align*}
\beta & = -({\varphi }^{\ast}_{\infty}\beta - {\varphi }^{\ast }_0\beta ) = 
- \int^{\infty}_0 \frac{\dee }{\dee t} ({\varphi }^{\ast }_t\beta ) \, \dee t \\
& = -\int^{\infty}_0 {\varphi }^{\ast }_t( L_X \beta ) \, \dee t = 
-\int^{\infty}_0 {\varphi }^{\ast }_t \big( \dee \, (X \lefthook \beta ) + 
X \lefthook \dee \beta \big) \, \dee t \\
& = -\int^{\infty}_0 {\varphi }^{\ast }_t \big( \dee \, (X \lefthook \beta )\big) \, \dee t, 
\, \, \mbox{since $\beta $ is closed} \\
& = -\dee \big( \int^{\infty}_0 {\varphi }^{\ast}_t(X \lefthook \beta ) \, \dee t \big), \, \, 
\mbox{since $\dee {\varphi }^{\ast }_t = {\varphi }^{\ast }_t \dee $\, .}
\end{align*}
The $(\ell -1)$-form $\alpha = 
\int^{\infty}_0 {\varphi }^{\ast}_t(X \lefthook \beta ) \, \dee t $ on $B$ is $H$ invariant, 
since ${\varphi }_t$ commutes with the $H$ action on $B$, and 
$X\lefthook \beta $ is an $H$ invariant $(\ell -1)$-form on $B$, because the vector field $X$ and the $\ell $-form $\beta $ are both $H$ invariant. Thus $\beta = \dee \alpha $ on $B$. Moreover, $\alpha $ is $G$ semi-basic, since 
\begin{displaymath}
L_{X_{\xi}}\big( {\varphi }^{\ast }_t(X \lefthook \beta ) \big) = 
{\varphi }^{\ast }_t \big( L_{X_{\xi }} (X \lefthook \beta ) \big) = 
{\varphi }^{\ast }_t \big( L_{X_{\xi }}X \lefthook L_{X_{\xi }}\beta \big) = 
0. 
\end{displaymath}
The last equality above follows because the $\ell $-form $\beta $ is $G$ semi-basic. \hfill $\square $ \medskip 

Since $M/G$ is a locally contractible space, we have \medskip 

\noindent \textbf{Proposition 14.} \textbf{(Poincar\'{e} lemma.)} \textit{Let 
$\overline{\theta}$ be a closed $\ell $-form on $M/G$ with $\ell \ge 1$. For each 
$\overline{m} \in M/G$ there is a contractible open neighborhood 
${\mathcal{U}}_{\, \overline{m}}$ of $\overline{m}$ 
and an $(\ell -1)$-form $\phi $ on ${\mathcal{U}}_{\, \overline{m}}$ such that 
$\overline{\theta}  = \dee \overline{\phi }$ on ${\mathcal{U}}_{\, \overline{m}}$.}  \medskip 

\noindent \textbf{Proof.} Since the $G$ action $\Phi $ on $M$ is proper, it has a slice $S_m$ at $m$, where $\pi (m) = \overline{m}$. Using Bochner's lemma there is an open neighborhood $U_m$ of $m$ in $S_m$, which is the image of an $H = G_m$ invariant open ball $B \subseteq T_mM$, centered at the origin $0_m$ whose closure is compact, under a diffeomorphism $\psi : B \subseteq T_mM \rightarrow U_m \subseteq 
M$. The diffeomorphism $\psi $ intertwines the linear $H$ action 
\begin{displaymath}
H \times T_m M \rightarrow T_m M:(h, v_m) \mapsto  T_m{\Phi }_gv_m
\end{displaymath}
with the $H \subseteq G$ action $\Phi $ on $U_m$. Let $\vartheta $ be 
the semi-basic $G$ invariant form on $G \cdot U_m$ such that 
$({\pi }_{|G \cdot U_m})^{\ast }\overline{\theta } = \vartheta $. Since $\overline{\theta }$ is closed by hypothesis, it follows that the semi-basic $\ell $-form $\vartheta $ on $G\cdot U_m$ is closed. Let $\phi = {\vartheta }|U_m$. Then $\phi $ is a semi-basic $H$ invariant closed $\ell $-form on $U_m$. Under the map $\psi $ the $\ell $-form $\phi $ pulls back to a 
$G$ semi-basic $H$ invariant $\ell $-form ${\psi }^{\ast }\phi $ on $B \subseteq T_mM$. By Lemma 12 there is a $G$ semi-basic $H$ invariant $(\ell -1)$-form 
$\gamma $ on $B$ such that 
${\psi }^{\ast }\phi = \dee \gamma $. Hence $\alpha = {\psi }_{\ast }\gamma $ 
is a semi-basic $H$ invariant $(\ell -1)$-form on $U_m$. The $(\ell -1)$-form 
$\alpha $ on $U_m$ extends to a $G$ invariant $(\ell -1)$-form 
$\sigma $ on $G \cdot U_m$ defined by 
\begin{displaymath}
\sigma ({\Phi }_g(m))(T_s{\Phi }_g v_s) = \alpha (s) v_s, 
\end{displaymath}
for every $s \in U_m$ and every $v_s \in T_sS_m$. Arguing as in the proof of 
Lemma 6, it follows that $\sigma $ is a smooth $G$ invariant 
$(\ell -1)$-form on $G \cdot U_m$. The form $\sigma $ is semi-basic. Moreover, $\dee \delta = \vartheta $ on $G \cdot U_m$, since for every $g \in G$ one has 
\begin{displaymath}
\dee \sigma = \dee \, ({\Phi }^{\ast }_g \alpha ) = {\Phi }^{\ast }_g(\dee \alpha ) = 
{\Phi }^{\ast }_g(\vartheta )  = \vartheta . 
\end{displaymath}
Let ${\mathcal{U}}_{\, \overline{m}} = \pi (U_m)$. Since $U_m$ is contractible and 
the $G$ orbit map $\pi $ is continuous and open, it follows that the open neighborhood ${\mathcal{U}}_{\, \overline{m}}$ of $\overline{m} \in M/G$ is contractible. Since the $\ell $-form $\sigma $ is semi-basic, there is an $\ell $-form $\overline{\phi}$ on ${\mathcal{U}}_{\, \overline{m}}$ such that 
${\pi }^{\ast }\overline{\phi } = \sigma $ on $G\cdot m$. On $G\cdot U_m$ we have 
\begin{displaymath}
{\pi }^{\ast }\overline{\theta } = \vartheta = \dee \sigma =  
\dee \, ({\pi }^{\ast }\overline{\phi }) = {\pi }^{\ast }(\dee \overline{\phi }). 
\end{displaymath}
Because the orbit map $\pi $ is surjective, it follows that $\overline{\theta } = \dee \overline{\phi }$ on ${\mathcal{U}}_{\, \overline{m}}$, which proves the proposition. \hfill $\square $ \medskip 

\noindent \textbf{Lemma 13.} \textit{Let $\overline{f} \in C^{\infty}(M/G)$ and 
suppose that ${\mathcal{U}}_{\, \overline{m}}$ is a connected open neighborhood 
of $\overline{m} \in M/G$ such that $\dee \overline{f} = 0$, then 
$\overline{f}$ is constant on ${\mathcal{U}}_{\, \overline{m}}$.} \medskip 

\noindent \textbf{Proof}. It follows from our hypotheses that 
$f = {\pi }^{\ast }\overline{f}$ is a smooth $G$ invariant function on 
the open connected component $U_m$ of 
${\pi }^{-1}({\mathcal{U}}_{\, \overline{m}})$ 
containing $m$. Moreover, on $U_m$ we have 
\begin{displaymath}
\dee f = \dee \, ({\pi }^{\ast }\overline{f} ) = {\pi }^{\ast }(\dee \overline{f}) =0.
\end{displaymath}
Since $M$ is a smooth manifold, it follows that $f$ is constant on $U_m$. Hence 
$\overline{f}$ is constant on the connected open set 
$\pi (U_m) ={\mathcal{U}}_{\, \overline{m}}$ because $\pi $ is a continuous open map. \hfill $\square $ \medskip 

To prove de Rham's theorem, we will need some sheaf theory, which 
can be found in appendix C of Lukina, Takens, and Broer 
\cite{lukina-takens-broer}. Let $\mathcal{U} = 
{ \{ {\mathcal{U}}_{\alpha } \} }_{\alpha \in I}$ be an open covering of 
$M/G$. Because $M/G$ is locally contractible, 
the open covering $\mathcal{U}$ has a good refinement ${\mathcal{U}}'$, 
that is, every ${\mathcal{U}}_{\beta } \in {\mathcal{U}}'$ with $\beta \in J 
\subseteq I$ is locally contractible and ${\mathcal{U}}_{{\beta }_1} \cap \cdots \cap {\mathcal{U}}_{{\beta }_n}$ is either contractible or empty. In addition, because $M/G$ is paracompact, every open covering has a locally finite subcovering. Since the $G$ action on $M$ is proper, the orbit space $M/G$ has 
a $C^{\infty}(M/G)$ partition of unity subbordinate to the covering $\mathcal{U}$. \medskip 

Define the differential exterior algebra valued sheaf $\Lambda $ over $M/G$ by 
\begin{displaymath}
\Lambda : {\mathcal{U}}_{\alpha } \mapsto \big( \Lambda (
{\mathcal{U}}_{\alpha }), \wedge , {\dee }_{\, | {\mathcal{U}}_{\alpha}} \big) , 
\end{displaymath}
whose sections are differential forms on ${\mathcal{U}}_{\alpha }$. The 
sheaf $\Lambda $ induces the subsheaves 
\begin{displaymath}
{\Lambda }^{\ell }: {\mathcal{U}}_{\alpha } \rightarrow \big( {\Lambda }^{\ell }(M/G), 
\wedge , {\dee }_{\, {\mathcal{U}}_{\alpha }} \big) , 
\end{displaymath}
whose sections are differential $\ell $-forms on ${\mathcal{U}}_{\alpha }$. Note that 
\begin{displaymath}
{\Lambda } \rightarrow M/G: {\Lambda }_{\overline{m}} = 
\sum_{\ell }{\Lambda }^{\ell }_{\overline{m}} \rightarrow \overline{m} 
\vspace{-.15in}
\end{displaymath}
is a smooth vector bundle, as is 
\begin{displaymath}
{\Lambda }^{\ell } \rightarrow M/G: {\Lambda }^{\ell }_{\overline{m}} 
\rightarrow \overline{m}.
\end{displaymath}
Let $\mathcal{R}$ be the sheaf of locally constant $\R $-valued functions on 
$M/G$. The two exact sequence of sheaves \medskip 
\begin{displaymath}
0 \rightarrow \mathcal{R} \rightarrow \Lambda \rightarrow \cdots  \, \, 
\mathrm{and} \, \, 0 \rightarrow \mathcal{R} \rightarrow {\Lambda }^{\ell } \rightarrow 
\cdots  
\end{displaymath} 
are exact. \medskip 

We say that the sheaf $\Lambda $ is fine if for every open subset $\mathcal{U}$ of $M/G$, every smooth function $\overline{f}$ on $M/G$ and every smooth 
section $s: \mathcal{U}\subseteq M/G \rightarrow \Lambda (\mathcal{U})$ of the sheaf $\Lambda $, then ${\overline{f}}_{|\mathcal{U}}\, \sigma \in \Lambda (\mathcal{U})$.  \medskip 

\noindent \textbf{Theorem 3.} \textit{The sheaves $\Lambda $ and 
${\Lambda }^{\ell }$ of sections of the vector bundles $\Lambda $ and 
${\Lambda }^{\ell }$ are fine.} \medskip 

\noindent \textbf{Proof.} We treat the case of the sheaf $\Lambda $. The 
proof for the sheaf ${\Lambda }^{\ell }$ is similar and is omitted. The definition of fineness holds by definition of differential form. \hfill $\square $ \medskip 

\noindent \textbf{Corollary 3} \textit{$\Lambda $ and ${\Lambda }^{\ell }$ are fine sheaves of sections over $M/G$, which is paracompact. Let 
$\mathcal{U}$ be an open covering of $M/G$. Then 
${\mathrm{H}}^q(\mathcal{U}, \Lambda )$, the sheaf of $q^{\mathrm{th}}$ cohomology group of $\mathcal{U}$ with 
values in the sheaf $\Lambda $, vanishes for all $q \in {\Z }_{\ge 1}$. Similarly, 
${\mathrm{H}}^q(\mathcal{U}, {\Lambda }^{\ell }) = 0$ for all $q \in {\Z}_{\ge 1}$.} 
\medskip 

We are now in position to formulate de Rham's theorem. 
Let ${\Lambda }^{\ell }$ be the sheaf of differential $\ell $-forms on $M/G$ and 
let $\dee : {\Lambda }^{\ell } \rightarrow {\Lambda }^{\ell +1}$ be the sheaf 
homomorphism induced by exterior differentiation. For each $\ell \in {\Z }_{\ge 0}$ 
let ${\mathcal{Z}}^{\ell } = \ker \dee \, $, whose elements are closed $\ell $-forms on $M/G$. By Lemma 13 ${\mathcal{Z}}^0 = \mathcal{R}$. Define the 
${\ell }^{\mathrm{th}}$ 
de Rham cohomology group ${\mathrm{H}}^{\ell }_{\mathrm{dR}}(M/G) = 
\Gamma (M/G, {\mathcal{Z}}^{\ell })\raisebox{-2pt}{\Large \mbox{$/$}} 
\dee \Gamma (M/G, {\Lambda }^{\ell -1})$ when $\ell \in {\Z }_{\ge 1}$ and 
${\mathrm{H}}^0_{\mathrm{dR}}(M/G) = \Gamma (M/G, {\mathcal{Z}}^0)$. Here 
$\Gamma (M/G, \mathcal{G})$ is the set of sections of the sheaf $M/G \rightarrow \mathcal{G}$. \medskip 

\noindent \textbf{Theorem 4} (\textbf{de Rham's theorem.}) \textit{The sheaf 
cohomology of ${\Lambda }^{\ell }$ with coefficients in $\mathcal{R}$ does not depend on the good covering $\mathcal{U}$ of $M/G$. Thus for every 
$\ell \in {\Z }_{\ge 0}$ the ${\ell }^{\mathrm{th}}$ de Rham cohomology group 
${\mathrm{H}}^{\ell }_{\mathrm{dR}}(M/G)$ is isomorphic to the 
${\ell }^{\mathrm{th}}$ sheaf cohomology group 
${\mathrm{H}}^{\ell }(\mathcal{U}, \mathcal{R})$ of the good covering 
$\mathcal{U}$ with values in the sheaf $\mathcal{R}$ of locally constant real valued functions.} \medskip 

\noindent \textbf{Proof.} We give a sketch, leaving out the homological algebra, 
which is standard. For more details, see \cite{lukina-takens-broer} or 
\cite{cushman-sniatycki21a}.  Let $\mathcal{U}$ be a good covering of 
$M/G$. The Poincar\'{e} lemma holds on any finite intersection of contractible open sets in $\mathcal{U}$, so the following sequence of sheaves is exact  
\begin{displaymath}
0 \rightarrow {\mathcal{Z}}^{\ell } \stackrel{\iota}{\longrightarrow} {\Lambda }^{\ell } 
\stackrel{\dee }{\longrightarrow} {\mathcal{Z}}^{\ell +1} \rightarrow 0,  
\end{displaymath}
where $\iota : {\mathcal{Z}}^{\ell } \rightarrow {\Lambda }^{\ell }$ is the inclusion mapping. This exact sequence gives rise to the long exact sequence of cohomology groups 
\begin{displaymath}
0 \longrightarrow {\mathrm{H}}^0(\mathcal{U}, {\mathcal{Z}}^{\ell }) 
\stackrel{{\iota}_{\ast }}{\longrightarrow} {\mathrm{H}}^0(\mathcal{U}, 
{\Lambda }^{\ell }) \stackrel{{\dee \, }_{\ast }}{\longrightarrow} 
{\mathrm{H}}^0(\mathcal{U}, {\mathcal{Z}}^{\ell +1}) 
\stackrel{{\delta }_{\ast }}{\longrightarrow } {\mathrm{H}}^1(\mathcal{U}, {\mathcal{Z}}^{\ell }) 
\longrightarrow \cdots , 
\end{displaymath}
where ${\iota }_{\ast }$, ${\dee \, }_{\ast }$, and ${\delta }_{\ast }$ are homomorphisms on cohomology induced by the \linebreak 
inclusion, exterior differentiation and coboundary 
homomorphisms, respectively. Since ${\Lambda }^{\ell}$ is a fine sheaf, its cohomology vanishes for $\ell \ge 1$ and the above sequence falls apart into the exact sequence 
\begin{equation}
0 \longrightarrow {\mathrm{H}}^0(\mathcal{U}, {\mathcal{Z}}^{\ell }) 
\stackrel{{\iota}_{\ast }}{\longrightarrow} {\mathrm{H}}^0(\mathcal{U}, 
{\Lambda }^{\ell }) 
\stackrel{{\dee \, }_{\ast }}{\longrightarrow} 
{\mathrm{H}}^0(\mathcal{U}, {\mathcal{Z}}^{\ell +1}) 
\stackrel{{\delta }_{\ast }}{\longrightarrow } {\mathrm{H}}^1(\mathcal{U}, {\mathcal{Z}}^{\ell }) 
\longrightarrow 0 , 
\label{eq-s5two}
\end{equation}
and for every $k \ge 1$ the exact sequence 
\begin{equation}
0 \longrightarrow {\mathrm{H}}^k(\mathcal{U}, {\mathcal{Z}}^{\ell +1}) 
\stackrel{{\delta }_{\ast }}{\longrightarrow} {\mathrm{H}}^{k+1}(\mathcal{U}, 
{\mathcal{Z}}^{\ell }) \longrightarrow 0.
\label{eq-s5three}
\end{equation}
Now ${\mathrm{H}}^0(M/G, \mathcal{R}) = \Gamma (M/G, {\mathcal{Z}}^0) = 
{\mathrm{H}}^0_{\mathrm{dR}}(M/G)$. Applying the sequence (\ref{eq-s5two}) 
consecutively gives 
\begin{displaymath}
{\mathrm{H}}^{\ell }(\mathcal{U}, \mathcal{R}) \simeq  
{\mathrm{H}}^1(\mathcal{U}, {\mathcal{Z}}^{\ell -1}). 
\end{displaymath}
Exactness of the sequence (\ref{eq-s5three}) gives 
\begin{displaymath}
{\mathrm{H}}^1(\mathcal{U}, {\mathcal{Z}}^{\ell -1}) 
\simeq {\mathrm{H}}^0(\mathcal{U}, {\mathcal{Z}}^{\ell } ) 
\raisebox{-2pt}{\Large \mbox{$/$}} 
\raisebox{-3pt}{${\dee }_{\, \ast }({\mathrm{H}}^0(\mathcal{U}, 
{\Lambda }^{\ell-1}))$} . 
\end{displaymath}
Here $\simeq $ means is isomorphic to. \hfill $\square $ \medskip 

\noindent \textbf{Corollary 4.} \textit{For the zeroth cohomology we have 
\begin{align}
&{\mathrm{H}}^{\ell }(\mathcal{U}, \mathcal{R}) \simeq  
\Gamma (M/G, {\mathcal{Z}}^{\ell }) \raisebox{-2pt}{\Large \mbox{$/$}} \hspace{-2pt}\raisebox{-2pt}{$\dee \, (\Gamma (M/G, {\Lambda }^{\ell -1}))$} = 
{\mathrm{H}}^{\ell }_{\mathrm{dR}}(M/G) \, \, \mbox{for $\ell \in {\Z}_{\ge 1}$.}
\tag*{$\square $}
\end{align} }
 
Our version of de Rham's theorem is not the same as 
Smith's version, since the ${\Z }_2$ invariant semi-basic $1$-form 
$x_2\dee x_1 - x_1 \dee x_2$ in section 6  is not a Smith $1$-form, 
see also Smith \cite[p.133]{smith}. However, his cohomology and ours agree. Our results extend those of Koszul \cite{koszul}, who hypothesized that $M/G$ was a smooth manifold and that the group $G$ was compact.  

\section{An example}

In this section we give an example, which illustrates Theorem 2 and the construction of differential $1$-forms on the orbit space of a proper group action on a smooth manifold. \medskip 

Consider the ${\Z }_2$ action on ${\R }^2$ generated by 
\begin{displaymath}
\zeta : {\R }^2 \rightarrow {\R }^2: x = (x_1,x_2) \mapsto (-x_1,-x_2) = -x. 
\end{displaymath}
The algebra of ${\Z }_2$ invariant polynomials on ${\R }^2$ is generated by 
the polynomials ${\sigma }_1(x) = x^2_1$, 
${\sigma }_2(x) = x^2_2$, and ${\sigma }_3(x) =x_1x_2 $, which are subject to the relation 
\begin{equation}
{\sigma }^2_3(x) = {\sigma }_1(x){\sigma }_2(x), \, \, {\sigma }_1(x) \ge 0 \, \, \& \, \, 
{\sigma }_2(x) \ge 0, \, \, \mbox{for all $x \in {\R }^2$.} 
\label{eq-exone}
\end{equation}
Let 
\begin{equation}
\sigma : {\R }^2 \rightarrow \Sigma \subseteq {\R }^3: x \mapsto 
\big( {\sigma }_1(x), {\sigma }_2(x), {\sigma }_3(x) \big) 
\label{eq-exthreestarnw}
\end{equation}
be the Hilbert map of the ${\Z}_2$ action associated to the polynomial generators ${\sigma }_1(x)$, ${\sigma }_2(x)$, and ${\sigma }_3(x)$. The map $\sigma $ (\ref{eq-exthreestarnw}) is the orbit map of the ${\Z }_2$ action on ${\R }^2$.
The relation 
\begin{equation}
{\sigma }^2_3 = {\sigma }_1{\sigma }_2 \, \, \, {\sigma }_1 \ge 0 \, \, \& \, \, 
{\sigma }_2 \ge 0
\label{eq-sec6doublezero}
\end{equation}
defines the orbit space ${\R }^2/{\Z}_2$ as a closed 
semialgebraic subset $\Sigma $ of ${\R }^3$ with coordinates 
$({\sigma }_1, {\sigma }_2, {\sigma }_3)$. Geometrically $\Sigma $ is a cone in 
${\R }^3$ with vertex $(0,0,0)$.  \medskip 

Because ${\Z }_2$ is a compact Lie group, which acts linearly on 
${\R }^2$, Schwarz' theorem \cite{schwarz} implies that the space 
$C^{\infty}({\R }^2)^{{\Z }_2}$ of ${\Z }_2$ invariant smooth functions on 
${\R }^2$ is equal to $C^{\infty}(\Sigma )$, where $\overline{f }\in 
C^{\infty}(\Sigma )$ if and only if there is an $\overline{F} \in C^{\infty}({\R }^3)$ such that $\overline{f} = {\overline{F}}_{|\Sigma }$.  \medskip 

\noindent \textbf{Lemma 14.} \textit{Let $f \in C^{\infty}({\R }^2)$ satisfy $f(x_1,x_2) =-f(-x_1,-x_2)$ for every $(x_1,x_2) \in {\R }^2$. Then there are $f_1$, $f_2 \in C^{\infty}({\R }^2)^{{\Z }_2}$ such that $f(x_1,x_2) = x_1f_1(x_1,x_2)+x_2f_2(x_1,x_2)$ for every $(x_1,x_2) \in {\R }^2$.} \medskip 

\noindent \textbf{Proof.} Since $f(x_1,x_2) = -f(-x_1,-x_2)$, it follows that $f(0,0) =0$. Suppose that there is an integer $k \ge 1$ such that $D^jf(0,0) = (0,0)$ for $0 \le j \le k-1$ and $D^kf(0,0) \ne (0,0)$. Then by Taylor's theorem with integral remainder we have $f(x,y) = 
\sum^k_{\ell =0} g_{\ell }(x_1,x_2)x^{k-\ell }_1x^k_2$ for every $(x_1,x_2) \in {\R }^2$, where 
$g_{\ell } \in C^{\infty}({\R }^2)$ and $g_{\ell }(0,0) = 
\frac{{\partial }^k f}{\partial x^{k-\ell }_1 \partial x^{\ell }_2}(0,0)$ for $0 \le \ell \le k$. By hypothesis
\begin{displaymath}
\begin{array}{rl}
\sum^k_{\ell =0} g_{\ell }(x_1,x_2)x^{k-\ell }_1x^k_2 & = f(x_1,x_2) \\
\rule{0pt}{12pt}&\hspace{-.5in}= -f(-x_1,-x_2) = 
(-1)^{k+1} \sum^k_{\ell =0} g_{\ell }(-x_1,-x_2)x^{k-\ell }_1x^k_2 .
\end{array}
\end{displaymath}
So 
\begin{equation}
g_{\ell }(x_1,x_2) = (-1)^{k+1} g_{\ell }(-x_1,-x_2), \, \, \mbox{for all $0 \le \ell \le k$.} 
\label{eq-exfive}
\end{equation}
If $k$ is odd Equation (\ref{eq-exfive}) implies 
$g_{\ell } \in C^{\infty}({\R}^2)^{{\Z }_2}$ for $0 \le \ell \le k$. Consequently, 
\begin{displaymath}
f(x_1,x_2) = x_1\big( g_0(x_1,x_2)x^{k-1}_1 \big) + 
x_2\big( \sum^k_{\ell =1}g_{\ell }(x_1,x_2)x^{k-\ell }_1x^{\ell -1}_2 \big) , 
\end{displaymath}
which proves the lemma when $k$ is odd. When $k$ is even, Equation 
(\ref{eq-exfive}) reads $g_{\ell }(x_1,x_2) = - g_{\ell }(-x_1,-x_2)$ for $0 \le \ell \le k$, which implies $g_{\ell }(0,0) =0$ for $0 \le \ell \le k$. In particular, $D^kf(0,0) =0$, which contradicts our hypothesis. 
\par Now suppose that $f$ is flat at $(0,0)$, that is, $D^kf(0,0) =0$ for every 
$k >0$. Then $x_1$ and $x_2$ divide $f$, that is, $f_1 = f/(2x_1)$ and $f_2 = f/(2x_2)$ are smooth functions on ${\R }^2$. To see this note that $f_1$ and $f_2$ are smooth for all $(x_1,x_2) \ne (0,0)$. Since $f$ is flat at $(0,0)$, so are $f_1$ and $f_2$. Clearly $f(x_1,x_2) = x_1f_1(x_1,x_2) +x_2f_2(x_1,x_2)$. From 
\begin{displaymath}
f_1(-x_1,-x_2) = f(-x_1,-x_2)/(2(-x_1)) = f(x_1,x_2)/(2x_1) = f_1(x_1,x_2)
\end{displaymath}
it follows that $f_1 \in C^{\infty}({\R }^2)^{{\Z }_2}$. Similarly, 
$f_2 \in C^{\infty}({\R }^2)^{{\Z}_2}$. \hfill $\square $ \medskip 

\noindent \textbf{Proposition 15.} \textit{The $C^{\infty}({\R }^2)^{{\Z }_2}$ module 
$\mathfrak{X}({\R }^2)^{{\Z }_2}$ of ${\Z }_2$ invariant smooth vector fields on 
${\R }^2$ is generated by 
\begin{equation}
X_1 = x_1\frac{\partial }{\partial x_1}, \, \, X_2 = x_2 \frac{\partial }{\partial x_1}, \, \, 
X_3 = x_1 \frac{\partial }{\partial x_2}, \, \, \mathrm{and} \, \, X_4 = x_2 \frac{\partial }{\partial x_2}. 
\label{eq-exfour}
\end{equation} }

\noindent \textbf{Proof.} A smooth vector field $X$ on ${\R }^2$ may be written as 
$X(x_1,x_2) = f(x_1,x_2) \frac{\partial }{\partial x_1} $ $+ g(x_1,x_2) \frac{\partial }{\partial x_2}$, where $f$ and $g$ are smooth. 
$X \in \mathfrak{X}({\R }^2)^{{\Z }_2}$ if and only if 
\begin{align*}
f(x_1,x_2) \frac{\partial }{\partial x_1} + g(x_1,x_2) \frac{\partial }{\partial x_2} & = 
X(x_1,x_2) \\ 
&\hspace{-1in} = {\zeta }^{\ast }X (x_1,x_2) = -f(-x_1,-x_2) \frac{\partial }{\partial x_1} - 
g(-x_1,-x_2) \frac{\partial }{\partial x_2}, 
\end{align*}
that is, $f(x_1,x_2) = -f(-x_1,-x_2)$ and $g(x_1,x_2)= -g(-x_1,-x_2)$ for every 
$(x_1,x_2) \in {\R }^2$. 
Using Lemma 14 write $f(x_1,x_2) = x_1 g_1(x_1,x_2) + x_2 g_2(x_1,x_2)$ and 
$g(x_1,x_2) = x_1 h_1(x_1,x_2) +x_2h_2(x_1,x_2)$, where $g_1$, $g_2$, $h_1$, and $h_2 \in C^{\infty}({\R }^2)^{{\Z }_2}$. Hence for every $(x_1,x_2) \in {\R }^2$ we have 
\begin{align}
X(x_1,x_2) & = (x_1 g_1(x_1,x_2) +x_2g_2(x_1,x_2)) \frac{\partial }{\partial x_1} \notag \\
&\hspace{.75in} + (x_1h_1(x_1,x_2) +x_2h_2(x_1,x_2)) \frac{\partial }{\partial x_2} 
\notag \\
& = (g_1 X_1 +g_2 X_2 +h_1X_3 +h_2 X_4)(x_1,x_2), \notag 
\end{align}
where $g_1$, $g_2$, $h_1$, and $h_2 \in C^{\infty}({\R }^2)^{{\Z }_2}$. \hfill $\square $ \medskip 

\noindent \textbf{Lemma 15} \textit{The vector fields 
\begin{equation}
\begin{array}{l}
Y_1 = 2{\overline{\sigma }}_1 \frac{\partial }{\partial {\sigma }_1} + 
{\overline{\sigma }}_3 \frac{\partial }{\partial {\sigma }_3}, 
\, \, Y_2 = 2{\overline{\sigma }}_3 \frac{\partial }{\partial {\sigma }_1} + 
{\overline{\sigma }}_2 \frac{\partial }{\partial {\sigma }_3} \\
\rowspace Y_3 = 2{\overline{\sigma }}_3 \frac{\partial }{\partial {\sigma }_2} + 
{\overline{\sigma }}_1\frac{\partial }{\partial {\sigma }_3}, 
\, \, Y_4 = 2{\overline{\sigma }}_2 \frac{\partial }{\partial {\sigma }_2} + 
{\overline{\sigma }}_3 \frac{\partial }{\partial {\sigma }_3} 
\end{array}
\label{eq-exsixa}
\end{equation}
on $\Sigma \subseteq {\R }^3$, where ${\overline{\sigma }}_i = 
({\sigma }_i)_{|\Sigma }$ for $i=1,2,3$, are $\sigma $ related to the 
${\Z }_2$ invariant vector fields $X_i$ (\ref{eq-exfour}) for $i=1,2,3$.} \medskip 

\noindent \textbf{Proof.} The calculation 
\begin{displaymath}
\begin{array}{lclcl}
L_{X_1}{\sigma }_1 = x_1\frac{\partial x^2_1}{\partial x_1} = 2{\sigma }_1, & &
L_{X_1}{\sigma }_2 = x_1\frac{\partial x^2_2}{\partial x_1} =0, & & 
L_{X_1}{\sigma }_3 = x_1\frac{\partial (x_1x_2)}{\partial x_1} = {\sigma }_3 \\
\rule{0pt}{15pt}L_{X_2}{\sigma }_1 = x_2\frac{\partial x^2_1}{\partial x_1} = 2{\sigma }_3, & &
L_{X_2}{\sigma }_2 = x_2\frac{\partial x^2_2}{\partial x_1} =0, & & 
L_{X_2}{\sigma }_3 = x_2\frac{\partial (x_1x_2)}{\partial x_1} = {\sigma }_2 \\
\rule{0pt}{15pt}L_{X_3}{\sigma }_1 = x_1\frac{\partial x^2_1}{\partial x_2} = 0, & &
L_{X_3}{\sigma }_2 = x_1\frac{\partial x^2_2}{\partial x_2} =2{\sigma }_3, & & 
L_{X_3}{\sigma }_3 = x_1\frac{\partial (x_1x_2)}{\partial x_2} = {\sigma }_1 \\
\rule{0pt}{15pt}L_{X_4}{\sigma }_1 = x_2\frac{\partial x^2_1}{\partial x_2} = 0, & &
L_{X_4}{\sigma }_2 = x_2\frac{\partial x^2_2}{\partial x_2} =2{\sigma }_2, & & 
L_{X_4}{\sigma }_3 = x_2\frac{\partial (x_1x_2)}{\partial x_2} = {\sigma }_3
\end{array}
\end{displaymath}
gives the vector fields 
\begin{equation}
\begin{array}{l}
{\widetilde{Y}}_1({\sigma }_1,{\sigma }_2, {\sigma }_3) = 2{\sigma }_1 \frac{\partial }{\partial {\sigma }_1} + 
{\sigma }_3 \frac{\partial }{\partial {\sigma }_3}, 
\, \, {\widetilde{Y}}_2({\sigma }_1,{\sigma }_2, {\sigma }_3) = 2{\sigma }_3 \frac{\partial }{\partial {\sigma }_1} + 
{\sigma }_2 \frac{\partial }{\partial {\sigma }_3} \\
\rowspace {\widetilde{Y}}_3({\sigma }_1,{\sigma }_2, {\sigma }_3) = 2{\sigma }_3 \frac{\partial }{\partial {\sigma }_2} + 
{\sigma }_1\frac{\partial }{\partial {\sigma }_3}, 
\, \, {\widetilde{Y}}_4({\sigma }_1,{\sigma }_2, {\sigma }_3) = 2{\sigma }_2 \frac{\partial }{\partial {\sigma }_2} + 
{\sigma }_3 \frac{\partial }{\partial {\sigma }_3} 
\end{array}
\label{eq-exsixb}
\end{equation}
on ${\R }^3$. Since $L_{X_i}({\sigma }^2_3 - {\sigma }_1{\sigma }_2) = 
0$ for $i=1,2,3,4$, the vector fields ${\widetilde{Y}}_i$ on ${\R }^3$ given by (\ref{eq-exsixb}) leave invariant the ideal $I$ of $C^{\infty}({\R }^3)$ generated by ${\sigma }^2_3 - {\sigma }_1{\sigma }_2$. Hence for each $i=1,2,3,4$ the 
vector field ${\widetilde{Y}}_i$ define the vector field $Y_i = 
{\widetilde{Y}}_i|\Sigma $ on $\Sigma $, which is given 
in Equation (\ref{eq-exsixa}). The vector fields $Y_i$ are $\sigma $ related to the 
${\Z }_2$ invariant vector fields $X_i$ (\ref{eq-exfour}) for $i=1,2,3,4$, because 
\begin{align}
Y_i\big( \sigma (x) \big) = ({\widetilde{Y}}_i|\Sigma )({\sigma }_1, {\sigma }_2, {\sigma }_3) = T_x \sigma \, X_i(x). \tag*{$\square $}
\end{align} 

Since the tangent to the Hilbert mapping $\sigma $ (\ref{eq-exthreestarnw}) is defined and is surjective, the tangent bundle $T\Sigma$ of the semialgebraic variety $\Sigma$ (\ref{eq-exone}) is the semialgebraic subset of ${\R }^7$ with coordinates $\big( {\sigma }_1, {\sigma}_2, {\sigma }_3, {\widetilde{Y}}_1, {\widetilde{Y}}_2,{\widetilde{Y}}_3, {\widetilde{Y}}_4 \big) $ defined by Equations (\ref{eq-exone}) and 
\begin{displaymath}
{\sigma }_3{\widetilde{Y}}_1 -{\sigma }_1{\widetilde{Y}}_2 = 0 \, \, \, \mathrm{and} \, \, \, 
{\sigma }_3{\widetilde{Y}}_3 - {\sigma }_1{\widetilde{Y}}_4 = 0.
\end{displaymath}

By Theorem 2 every smooth vector field 
on $\Sigma$ is $\sigma $ related to a smooth ${\Z }_2$ invariant vector field on 
${\R }^2$. Because the $C^{\infty}({\R }^2)^{{\Z}_2}$ module 
$\mathfrak{X}({\R }^2)^{{\Z}_2}$ of smooth ${\Z}_2$ 
invariant vector fields on ${\R }^2$ is generated by the vector fields $X_i$ for $1 \le i \le 4$ given by Equation (\ref{eq-exfour}), it follows that the $\sigma $ related vector fields $Y_i$ for $1 \le i \le 4$ given by Equation (\ref{eq-exsixa}) generate the $C^{\infty}(\Sigma)$ module $\mathfrak{X}(\Sigma)$ of smooth vector fields on $\Sigma $. \medskip 

\noindent \textbf{Lemma 16.} \textit{The differential $1$-forms 
\begin{equation}
{\widetilde{\vartheta }}_1 = x_1 \dee x_1 , \, \, {\widetilde{\vartheta }}_2 = x_1 \dee x_2, \, \, 
{\widetilde{\vartheta }}_3 = x_2 \dee x_1 , \, \, {\widetilde{\vartheta }}_4 = x_2 \dee x_2 . 
\label{eq-s4twostarnw}
\end{equation}
generate the $C^{\infty}({\R }^2)^{{\Z }_2}$ module 
${\Lambda }^1({\R }^2)^{{\Z }_2}$ of ${\Z }_2$ invariant $1$-forms on ${\R }^2$.}  \medskip 

\noindent \textbf{Proof.} We 
use the differential forms 
\begin{equation}
\begin{array}{lcl}
{\vartheta }_1 = \dee x^2_1 = 2x_1 \dee x_1, & \quad &  
{\vartheta }_2 = \dee x^2_2 = 2x_2  \dee x_2, \\
\rule{0pt}{13pt}{\vartheta }_3 = \dee \,  (x_1x_2) = x_1 \dee x_2 + x_2 \dee x_1, & \quad  &
{\vartheta }_4 = x_1 \dee x_2 - x_2 \dee x_1 
\end{array} 
\label{eq-s4twostardaggernw}
\end{equation}
instead of those given in (\ref{eq-s4twostarnw}), because we then get 
${\vartheta }_k = 
\dee {\sigma }_k$ for $k=1,2,3$. Suppose that the $1$-form $\vartheta (x_1,x_2) = f_1(x_1,x_2) \dee x_1 $ $ + f_2(x_1,x_2) \dee x_2$ on ${\R}^2$, 
where $f_i \in C^{\infty}({\R }^2)$ for 
$i =1,2$, is invariant under the ${\Z }_2$ action generated by 
$\zeta : {\R }^2 \rightarrow {\R }^2:(x_1,x_2) \mapsto (-x_1,-x_2)$. Then for every 
$(x_1,x_2) \in {\R }^2$  
\begin{align*}
f_1(x_1,x_2) \dee x_1 + f_2(x_1,x_2) \dee x_2 & = 
\vartheta (x_1,x_2 ) = ({\zeta }^{\ast }\vartheta )(x_1,x_2) \\
&\hspace{-1in} = f_1(-x_1,-x_2) \dee \, (-x_1) +f_2(-x_1,-x_2) \dee \, (-x_2) \\
&\hspace{-1in} = -f_1(-x_1,-x_2) \dee x_1 - f_2(-x_1,-x_2) \dee x_2. 
\end{align*}
So ${\zeta }^{\ast }\vartheta = \vartheta $ if and only if for $i=1,2$ one has 
$f_i(x_1,x_2) = -f_i(-x_1,-x_2)$ for every $(x_1,x_2) \in {\R }^2$. 
By Lemma 14 if  $g(x_1,x_2) = -g(-x_1,-x_2)$ for some $g \in C^{\infty}({\R }^2)$, then there are $g_i \in C^{\infty}({\R }^2)^{{\Z}_2}$ for $i=1,2$ such that 
$g(x_1,x_2) = x_1g_1(x_1,x_2) + x_2 g_2(x_1,x_2)$ for every 
$(x_1,x_2) \in {\R }^2$. Consequently, for some $h_1$, $h_2$, $k_1$, $k_2 \in 
C^{\infty}({\R }^2)^{{\Z}_2}$ 
\begin{align*}
\vartheta (x_1,x_2) & = \big( x_1h_1(x_1,x_2)  +x_2h_2(x_1,x_2) \big) \dee x_1 + 
\big( x_1k_1(x_1,x_2)  +x_2k_2(x_1,x_2) \big) \dee x_2 \\
& = h_1(x_1,x_2)\, x_1 \dee x_1 + h_2(x_1,x_2)\, x_2 \dee x_1  \\
& \hspace{.5in} + k_1(x_1,x_2)\, x_1 \dee x_2 + k_2(x_1,x_2)\, x_2 \dee x_2 \\
& =  {\widetilde{h}}_1(x_1,x_2)\, {\vartheta }_1 + {\widetilde{h}}_2(x_1,x_2)\, {\vartheta}_2 + {\widetilde{k}}_1(x_1,x_2)\, {\vartheta }_3 + 
{\widetilde{k}}_2(x_1,x_2)\, {\vartheta }_4, 
\end{align*}
for every $(x_1,x_2) \in {\R }^2$. Here ${\widetilde{h}}_1 = 2h_1$, 
${\widetilde{h}}_2 = 2k_2$, ${\widetilde{k}}_1 = k_1+h_2$, and 
${\widetilde{k}}_2 =  k_1 - h_2$. This proves the lemma. \hfill $\square $ \medskip 

For $i=1, \ldots , 4$ define the $1$-forms ${\theta }_i$ on $\Sigma $ by 
\begin{equation}
{\sigma }^{\ast }\big( (Y_i |\, \Sigma )  \lefthook {\theta }_i \big) = X_i \lefthook {\vartheta }_i, 
\label{eq-s4doublezero}
\end{equation}
see the proof of Proposition 10.  The $1$-forms ${\theta }_i$  generate the 
$C^{\infty}({\R }^2/{\Z}_2)$ module of $1$-forms on $\Sigma $, since the 
${\Z}_2$ invariant $1$-forms ${\vartheta }_i$ for $i=1, \ldots , 4$ generate the 
$C^{\infty}({\R }^2)^{{\Z }_2}$ module ${\Lambda }({\R }^2)^{{\Z }_2}$ 
of ${\Z }_2$ invariant $1$-forms on ${\R }^2$, 
see Lemma 14. Every ${\Z }_2$ invariant $1$-form on ${\R }^2$ is semi-basic, since the Lie algebra of ${\Z }_2$ is $\{ 0 \}$. \medskip  

\noindent \textbf{Fact 1.} \textit{On $\Sigma $ we have 
\begin{align}
 {\theta }_1 & = \dee {\overline{\sigma}}_1, \, \, {\theta }_2 = \dee {\overline{\sigma }}_2, \, \, 
 \mathrm{and} \, \, {\theta }_3 = \dee {\overline{\sigma }}_3 . 
 \label{eq-s4nwzero} 
 \end{align}
Let ${\theta }_4$ be the $1$-form on $\Sigma $ defined by its values 
\begin{equation}
\begin{array}{rl}
Y_1|\Sigma \lefthook {\theta }_4 = -{\overline{\sigma }}_3, & 
Y_2|\Sigma \lefthook {\theta }_4 = {\overline{\sigma }}_1, \\ 
\rule{0pt}{12pt}Y_3|\Sigma \lefthook {\theta }_4 = -{\overline{\sigma }}_2, & 
Y_4|\Sigma \lefthook {\theta }_4 = -{\overline{\sigma }}_3. 
\end{array} 
\label{eq-s4zerostarnw}
\end{equation}
Here ${\overline{\sigma }}_i = {{\sigma }_i}_{|\Sigma }$ for $i=1,2,3$. 
The $1$-form ${\theta }_4$ is not the restriction of a $1$-form on 
${\R }^3$ to $\Sigma $.}  \medskip 

\noindent \textbf{Proof.} Equation (\ref{eq-s4nwzero}) follows immediately from the definition of ${\theta }_i$ given in Equation (\ref{eq-s4doublezero}). \medskip 

We give three proofs of the assertion about ${\theta }_4$. \medskip 

\noindent $\mathbf{1}$. Consider the $1$ form 
$\theta = \frac{{\sigma }_1}{2{\sigma }_3} \dee {\sigma }_2 - 
\frac{{\sigma }_2}{2{\sigma }_3} \dee {\sigma }_1$ on ${\R }^3$. Then 
\begin{align*}
{\sigma }^{\ast }({\theta }_{|\Sigma }) & = 
\frac{x^2_1 \dee x^2_2 - x^2_2 \dee x^2_1}{2 x_1x_2} 
= x_1 \dee x_2 - x_2 \dee x_1 = {\vartheta }_4.
\end{align*}
The following argument shows that the $1$ -form 
\begin{displaymath}
{\theta }_{|\Sigma } 
= \frac{{\sigma }_1}{2{\sigma }_3}_{|\Sigma } (\dee {\sigma }_2)_{|\Sigma } -
\frac{{\sigma }_2}{2{\sigma }_3}_{|\Sigma } (\dee {\sigma }_1)_{|\Sigma }
\end{displaymath} 
is not smooth, because its coefficients are not smooth functions on 
$\Sigma $. First we need some geometric information about the 
${\Z}_2$ orbit space $\Sigma \subseteq {\R }^3$ defined by 
${\sigma }^2_3 = {\sigma }_1{\sigma }_2$ with 
${\sigma }_1 \ge 0$ and ${\sigma }_2 \ge 0$. The only subgroups of 
${\Z }_2$ are the identity $\{ e \} $ and ${\Z }_2$. The isotropy group 
$({\Z }_2)_{x}$ at $x \in {\R }^2$ is ${\Z }_2$ if $x =0$ and 
$\{ e \} $ if $x \ne 0$. The corresponding orbit types are $\{ 0 \} $ and 
${\R }^2 \setminus \{ 0 \}$, whose image under the orbit map $\sigma $ 
is $\mathrm{O} = \{ (0,0,0) \} $, the vectex of the cone $\Sigma $, and 
$\Sigma \setminus \mathrm{O}$, which is a smooth manifold. Thus 
${\theta }_{|(\Sigma \setminus \mathrm{O})}$ is a smooth $1$-form, whose 
pull back under $\sigma $ is the smooth $1$-form ${\vartheta }_4$ on 
${\R }^2\setminus \{ (0,0) \}$. The $1$-form 
${\theta }_{|(\Sigma \setminus \mathrm{O})}$ does 
not extend to a smooth $1$-form ${\theta }_{|\Sigma }$ because 
the functions $\frac{{\sigma }_1}{2{\sigma }_3}_{|\Sigma }$ and 
$\frac{{\sigma }_2}{2{\sigma }_3}_{|\Sigma }$ are not smooth at 
$(0,0,0)$, the vertex of the cone $\Sigma $. To see this 
let ${\sigma }^0 = ({\sigma }^0_1, {\sigma}^0_2, {\sigma }^0_3) 
\in \Sigma \setminus \mathrm{O}$. The closed line segment 
${\ell }_{{\sigma }^0}([0,1])$, where 
\begin{displaymath}
{\ell }_{{\sigma }^0}: [0,1] \rightarrow \Sigma : t \mapsto t{\sigma }^0 = 
(t{\sigma }^0_1, t{\sigma}^0_2, t{\sigma }^0_3),
\end{displaymath}
lies in $\Sigma $ and joins $(0,0,0)$ to ${\sigma }^0$. Now 
$\frac{{\sigma }_1}{2{\sigma }_3}\big( {\ell }_{{\sigma }^0}(t) \big) =
\frac{t{\sigma }^0_1}{2t{\sigma }^0_3} = \frac{{\sigma }^0}{2{\sigma }^3}$. 
So $\frac{{\sigma }_1}{2{\sigma }_3}(0,0,0) = \frac{{\sigma }^0}{2{\sigma }^0_3}$. 
Hence the function $\frac{{\sigma }_1}{2{\sigma }_3}_{|\Sigma }$ is 
not continuous at $(0,0,0)$. A similiar argument shows that the function 
$\frac{{\sigma }_2}{2{\sigma }_3}_{|\Sigma }$ is not continuous 
at $(0,0,0)$. \medskip 

\noindent $\mathbf{2}$. The following argument shows that 
the $1$-form ${\theta}_4$ on $\Sigma $ defined in Equation 
(\ref{eq-s4zerostarnw}) is not the restriction to $\Sigma $ of 
any smooth $1$-form on ${\R }^3$. Suppose it is. Then 
${\theta }_4 = \sum^3_{j=1}{\overline{A}}_j \dee {\sigma }_j$,  
for some ${\overline{A}}_j \in C^{\infty}(\Sigma ) = C^{\infty}({\R }^3)/I$, 
where $I$ is the ideal of $C^{\infty}({\R }^3)$ generated by 
${\sigma }^2_3 - {\sigma }_1{\sigma }_2$. Using (\ref{eq-s4zerostarnw})  we get 
\begin{align*}
-{\sigma }_3 +I & = {\overline{\sigma }}_3 = 
Y_1|\Sigma \lefthook {\theta }_4 = ( 2{\sigma }_1\frac{\partial }{\partial {\sigma }_1} 
+ {\sigma }_3\frac{\partial }{\partial {\sigma }_3})|\Sigma \lefthook {\theta }_4 
= 2{\overline{\sigma }}_1 {\overline{A}}_1 + {\overline{\sigma }}_3 {\overline{A}}_3, \notag 
\end{align*}
which implies 
\begin{subequations}
\begin{equation}
-{\sigma }_3 = 2{\sigma }_1 A_1 + {\sigma}_3 A_3 + I. 
\label{eq-s4threeA}
\end{equation}
Similarly, 
\begin{align}
{\sigma }_1 & = 2{\sigma }_3A_1 +{\sigma }_2A_3 + I 
\label{eq-s4threeB} \\
-{\sigma }_2 & = 2{\sigma }_3 A_2 + {\sigma }_1 A_3 + I 
\label{eq-s4threeC} \\
-{\sigma }_3 & = 2{\sigma }_2A_2 + {\sigma }_3A_3 + I. 
\label{eq-s4threeD}
\end{align}
\end{subequations}
Set $A_1 = -{\sigma }_2$ and $A_3 = -1 + 2{\sigma }_3$. Then 
\begin{displaymath}
2{\sigma }_1A_1 + {\sigma }_3A_3 = -2{\sigma }_1{\sigma }_2 -{\sigma }_3+2{\sigma }^2_3 
= -{\sigma }_3 +I. 
\end{displaymath}
So Equation (\ref{eq-s4threeA}) holds. Multiplying (\ref{eq-s4threeB}) by 
${\sigma }_1$ and (\ref{eq-s4threeC}) by ${\sigma }_2$ and adding gives
\begin{align}
{\sigma }^2_1 - {\sigma }^2_2 & = 2({\sigma }_1A_1+{\sigma }_2A_2){\sigma }_3 + 2{\sigma }_1{\sigma }_2 A_3 +I \notag \\
& = 2{\sigma }_3[ {\sigma }_1A_1+{\sigma }_2A_2 +{\sigma }_3A_3] + I \notag \\
& = 2{\sigma }_3[ -{\sigma }_1{\sigma }_2 +{\sigma }_2A_2 - {\sigma }_3 + {\sigma }^2_3] + I 
\notag \\
& = 2{\sigma }_3[{\sigma }_2A_2 - {\sigma }_3 +{\sigma }^2_3] + I. 
\label{eq-s4threeE}
\end{align}
But ${\sigma }_3$ does not divide ${\sigma }^2_1-{\sigma }^2_2$, which does not lie in $I$. Thus Equation (\ref{eq-s4threeE}) does not hold for any choice of $A_2 \in C^{\infty}({\R }^3)$. Hence our hypothesis is false, that is, the $1$-form ${\theta }_4$ on ${\R }^2/{\Z }_2 = \Sigma $ is not the restriction to $\Sigma $ of a $1$-form on ${\R }^3$.  \medskip 

\noindent $\mathbf{3}$. Our third proof is more analytic. The 
$1$-form ${\theta }_4$ (\ref{eq-s4zerostarnw}) on 
the orbit space $\Sigma \subseteq {\R }^3$ is not the restriction to $\Sigma $ of a $1$ form $\theta = \sum^3_{j=1}A_j \dee {\sigma }_j$ on ${\R }^3$, where 
$A_j \in C^{\infty}({\R }^3)$.  Suppose 
that ${\theta }_4  = \theta |\Sigma $, then 
\begin{displaymath}
{\sigma  }^{\ast }(\dee {\theta }_4) = \dee \, ({\sigma }^{\ast }{\theta }_4) = \dee {\vartheta }_4 = \dee \, ( x_1 \dee x_2 - x_2 \dee x_1) = 2 \dee x_1 \wedge 
\dee x_2, 
\end{displaymath}
which does not vanish at $(0,0) \in {\R }^2$. However, the $2$-form 
\begin{align*}
{\sigma }^{\ast }( (\dee \theta ) |\Sigma ) & = 
\sum^3_{j=1} \dee \, ( {\sigma }^{\ast }A_j ) \wedge
{\sigma }^{\ast }(\dee {\overline{\sigma }}_j) \\
& =  \sum^3_{j=1} \dee \, ( {\sigma }^{\ast }A_j ) \wedge
\dee \, ({\sigma }^{\ast }({\overline{\sigma }}_j)) = 
\sum^3_{j=1}\dee \, ({\sigma }^{\ast }A_j) \wedge {\vartheta }_j
\end{align*}
vanishes at $(0,0)$, since the $1$-forms ${\vartheta }_j$ 
(\ref{eq-s4twostarnw}) for $j=1,2,3$ vanish at $(0,0)$. This is a contradiction, since $\dee {\theta }_4 = \dee \, (\theta |\Sigma )$.  \hfill $\square $ \bigskip  

\section{Acknowledgment}

The authors are grateful to Professor Gerald Schwarz for pointing out an \linebreak
error in the proof of Lemma 3.1, and for suggesting the inclusion in the 
bibliography of three additional papers \cite{schwarz80}, 
\cite{kankaanrinta}, and \cite{karshon-watts}, which are related to the problem under consideration.

\end{document}